# AN EFFECTIVE CRITERION AND A NEW EXAMPLE FOR BALLISTIC DIFFUSIONS IN RANDOM ENVIRONMENT

BY LAURENT GOERGEN

*ETH Zurich*

In the setting of multidimensional diffusions in random environment, we carry on the investigation of condition $(T')$, introduced by Sznitman [*Ann. Probab.* **29** (2001) 723–764] and by Schmitz [*Ann. Inst. H. Poincaré Probab. Statist.* **42** (2006) 683–714] respectively in the discrete and continuous setting, and which implies a law of large numbers with nonvanishing limiting velocity (ballistic behavior) as well as a central limit theorem. Specifically, we show that when $d \geq 2$, $(T')$ is equivalent to an effective condition that can be checked by local inspection of the environment. When $d = 1$, we prove that condition $(T')$ is merely equivalent to almost sure transience. As an application of the effective criterion, we show that when $d \geq 4$ a perturbation of Brownian motion by a random drift of size at most $\varepsilon > 0$ whose projection on some direction has expectation bigger than $\varepsilon^{2-\eta}, \eta > 0$, satisfies condition $(T')$ when $\varepsilon$ is small and hence exhibits ballistic behavior. This class of diffusions contains new examples of ballistic behavior which in particular do not fulfill the condition in [*Ann. Inst. H. Poincaré Probab. Statist.* **42** (2006) 683–714], (5.4) therein, related to Kalikow's condition.

**1. Introduction.** Diffusions in random environment emerged about 25 years ago from homogenization theory in the study of disordered media; see, for instance, [3]. Within the rich field of "random motions in random media," they are closely related to the discrete model of "random walks in random environment"; see [9, 22].

In the one-dimensional discrete setting a complete characterization of ballistic behavior, which refers to the situation where the motion tends to infinity in some direction with nonvanishing velocity, was established already in 1975 by Solomon [15]; see also [6, 10]. In the multidimensional setting, however, such a characterization has not been found yet, but a great deal of progress has been made over the last seven years: the so-called conditions $(T)$ and $(T')$ introduced by Sznitman (see [18, 19]) for random walks in random environment seem to be promising candidates for an equivalent description of ballistic behavior when the space dimension $d \geq 2$. In essence, one possible formulation of condition $(T)$ [see (1.12)] requires exponential decay of the probability that the trajectory exits a slab of growing width through one side rather than the other. These conditions have interesting consequences such as a ballistic law of large numbers and a central limit theorem. Their analogues in the setting of diffusions have been developed by Schmitz (see [11, 12]), and he used a previous result of Shen [13, 14] to show that they imply the same asymptotic behavior as mentioned before in the discrete setup. The drawback of the definitions of conditions $(T)$ or $(T')$ as they were stated in [11] [see (1.10)] is their asymptotic nature which makes them difficult to check by local considerations. To remedy this problem, we provide in the first part of this article an effective criterion, in the spirit of [19], which is equivalent to $(T')$ (see Theorem 2.6), and which can be checked by inspection of the environment in a finite box.

In the second part of this work, which is related to [20] in the discrete setting, we use the effective criterion to show that when $d \geq 4$, Brownian motion perturbed with a small random drift satisfying the assumption (1.16), fulfills condition $(T')$; see Theorem 3.1. As we will see below, this class of diffusions contains new examples for ballistic behavior beyond prior knowledge.

Before we discuss our results any further, we first describe the model. The random environment is specified by a probability space $(\Omega, \mathcal{A}, \mathbb{P})$ on which acts a jointly measurable group $\{t_x; x \in \mathbb{R}^d\}$ of $\mathbb{P}$-preserving transformations, with $d \geq 1$. The diffusion matrix and the drift of the diffusion in random environment are stationary functions $a(x, \omega), b(x, \omega), x \in \mathbb{R}^d, \omega \in \Omega$, with respective values in the space of nonnegative $d \times d$ matrices and in $\mathbb{R}^d$, that is,

$$
\begin{aligned}
a(x+y, \omega) &= a(x, t_y \omega), \\
b(x+y, \omega) &= b(x, t_y \omega) \qquad \text{for } x, y \in \mathbb{R}^d, \omega \in \Omega.
\end{aligned}
\tag{1.1}
$$

We assume that these functions are bounded and uniformly Lipschitz, that is, there is a $\bar{K} > 1$, such that for $x, y \in \mathbb{R}^d, \omega \in \Omega$,

$$
\begin{aligned}
|b(x, \omega)| + |a(x, \omega)| &\leq \bar{K}, \\
|b(x, \omega) - b(y, \omega)| + |a(x, \omega) - a(y, \omega)| &\leq \bar{K}|x - y|,
\end{aligned}
\tag{1.2}
$$

where $|\cdot|$ denotes the Euclidean norm for vectors and matrices. Further, we assume that the diffusion matrix is uniformly elliptic, that is, there is a $\nu > 1$ such that for all $x, y \in \mathbb{R}^d$, $\omega \in \Omega$:

$$
\frac{1}{\nu}|y|^2 \leq y \cdot a(x, \omega) y \leq \nu |y|^2. \tag{1.3}
$$

The coefficients $a, b$ satisfy a condition of finite range dependence: for $A \subset \mathbb{R}^d$, we define

(1.4) $$\mathcal{H}_A = \sigma(a(x, \cdot), b(x, \cdot); x \in A),$$

and assume that for some $R > 0$,

(1.5) $\mathcal{H}_A$ and $\mathcal{H}_B$ are independent under $\mathbb{P}$ whenever $d(A, B) \geq R$,

where $d(A, B)$ is the mutual Euclidean distance between $A$ and $B$. With the above regularity assumptions on $a$ and $b$, for any $\omega \in \Omega$, $x \in \mathbb{R}^d$, the martingale problem attached to $x$ and the operator

(1.6) $$\mathcal{L}_\omega = \tfrac{1}{2} \sum_{i,j=1}^d a_{ij}(\cdot, \omega) \partial_{ij}^2 + \sum_{i=1}^d b_i(\cdot, \omega) \partial_i$$

is well posed; see [17] or [2], page 130. The corresponding law $P_{x,\omega}$ on $C(\mathbb{R}_+, \mathbb{R}^d)$, unique solution of the above martingale problem, describes the diffusion in the environment $\omega$ and starting from $x$. We write $E_{x,\omega}$ for the expectation under $P_{x,\omega}$ and we denote the canonical process on $C(\mathbb{R}_+, \mathbb{R}^d)$ with $(X_t)_{t \geq 0}$. Observe that $P_{x,\omega}$ is the law of the solution of the stochastic differential equation

(1.7) $$\begin{aligned} dX_t &= \sigma(X_t, \omega) \, d\beta_t + b(X_t, \omega) \, dt, \\ X_0 &= x, \qquad P_{x,\omega}\text{-a.s.}, \end{aligned}$$

where, for instance, $\sigma(\cdot, \omega)$ is the square root of $a(\cdot, \omega)$ and $\beta$ is some $d$-dimensional Brownian motion under $P_{x,\omega}$. The laws $P_{x,\omega}$ are usually called "quenched laws" of the diffusion in random environment. To restore translation invariance, we consider the so-called "annealed laws" $P_x$, $x \in \mathbb{R}^d$, which are defined as semi-direct products:

(1.8) $$P_x \stackrel{\text{def}}{=} \mathbb{P} \times P_{x,\omega}.$$

Of course the Markov property is typically lost under the annealed laws.

We now come back to the object of this work. We start by recalling the definition of conditions $(T)$ and $(T')$ as stated in [11]. These conditions are expressed in terms of another condition $(T)_\gamma$ defined as follows. For a unit vector $\ell$ of $\mathbb{R}^d$, $d \geq 1$, and any $u \in \mathbb{R}$, consider the stopping times

(1.9) $$T_u^\ell = \inf\{t \geq 0; X_t \cdot \ell \geq u\}, \qquad \tilde{T}_u^\ell = \inf\{t \geq 0; X_t \cdot \ell \leq u\}.$$

For $\gamma \in (0, 1]$, we say that condition $(T)_\gamma$ holds relative to $\ell$, in shorthand notation $(T)_\gamma \mid \ell$, if for all unit vectors $\ell'$ in some neighborhood of $\ell$ and for all $b > 0$,

(1.10) $$\limsup_{L \to \infty} L^{-\gamma} \log P_0[\tilde{T}_{-bL}^{\ell'} < T_L^{\ell'}] < 0.$$

Condition $(T')$ relative to $\ell$ is then the requirement that

(1.11) $\qquad$ (1.10) holds for all $\gamma \in (0,1)$,

and condition $(T)$ relative to $\ell$ refers to the case where

(1.12) $\qquad$ (1.10) holds for $\gamma = 1$.

It is clear that $(T)$ implies $(T')$ and we show in Theorem 2.6 that $(T')$ is equivalent to $(T)_\gamma$ when $\gamma \in (\frac{1}{2}, 1)$. Moreover, it is conjectured that the conditions $(T)_\gamma$, $\gamma \in (0,1]$ are all equivalent.

Let us also mention that for $\gamma \in (0,1]$, condition $(T)_\gamma$ relative to $\ell$ is in essence equivalent to almost sure transience in direction $\ell$ together with finiteness of a stretched exponential moment of the size of the trajectory up to a certain regeneration time; see [11], Theorem 3.1 therein or [19] for a similar result in the discrete setting. The latter formulation of condition $(T')$ is especially appropriate to study the asymptotic properties of the diffusion. Indeed Schmitz showed in [11], Theorem 4.5 (see also [18]) that when $d \geq 2$, it enables us to verify the sufficient conditions of [13, 14] for a ballistic law of large numbers and a central limit theorem. However, the more geometrical expression (1.10) is better suited for our present purpose.

Despite the interest of the two above mentioned formulations of condition $(T')$, they are not "effective conditions" that can be checked by local inspection of the environment. Concrete examples where $(T')$ holds, besides the easy case where the projection of the drift on some unit vector is uniformly bounded away from 0 (see [11], Proposition 5.1), originate from a stronger condition going back to Kalikow; see [7, 21]. For instance, it is shown in [11], Theorem 5.2, and [12], Theorem 2.1 that there exists a constant $c_e > 0$ depending only on $\bar{K}, \nu, R, d$ [see (1.2)–(1.5)], such that condition $(T)$ holds when

(1.13) $\qquad \mathbb{E}[(b(0,\omega) \cdot \ell)_+] \geq c_e \mathbb{E}[(b(0,\omega) \cdot \ell)_-].$

In the first part of this work we derive an effective criterion in the above sense. We show (see Theorem 2.6) that when $d \geq 2$ for any direction $\ell$, $(T')|\ell$ is in essence equivalent to

(1.14) $\qquad \inf_{B, a \in (0,1]} \{c(d)\tilde{L}^{d-1} L^{3(d-1)+1} \mathbb{E}[\rho_B^a]\} < 1,$

with

(1.15) $\qquad \rho_B = \dfrac{P_{0,\omega}[X_{T_B} \notin \partial_+ B]}{P_{0,\omega}[X_{T_B} \in \partial_+ B]},$

provided in the above infimum, $B$ runs over all large boxes transversal to $\ell$ consisting of the points $x$ with $x \cdot \ell \in (-L + R + 2, L + 2)$ and other coordinates in an orthonormal basis with first vector $\ell$, smaller in absolute value than $\tilde{L}$, for $L \geq c'(d), R + 2 \leq \tilde{L} < L^3$. In the above formula for $\rho_B$, $T_B$ denotes the exit time from $B$ and $\partial_+ B$ is the part of the boundary of $B$ where $x \cdot \ell = L + 2$. The proof of Theorem 2.6 follows the strategy of Sznitman [19] and the sufficiency of the effective criterion is obtained by an induction argument along a growing sequence of boxes $B_k$ that tend to look like infinite slabs and in which suitable moments of $\rho_{B_k}$ are used to control moments of $\rho_{B_{k+1}}$. This allows us to deduce the asymptotic exit behavior (1.10) from slabs. As a first application of the effective criterion, we show the equivalence between $(T')$ and $(T)_\gamma$ when $\gamma \in (\frac{1}{2}, 1)$. Note also that $\rho_B$ in (1.15) reminds us of the decisive quantity appearing in the one-dimensional theorem of Solomon [15]. We will see in Section 2.1 that when $d = 1$, the box $B$ is replaced with an interval $(-L, L)$ and the existence of an $a \in (0, 1]$ and a $L > R$ such that $\mathbb{E}[\rho^a_{(-L,L)}] < 1$ is equivalent to $(T')$ and $(T)$ as well as to almost sure convergence to $+\infty$. Hence in opposition to the multidimensional case, condition $(T)$ does not imply ballistic behavior when $d = 1$.

In the second part of this article we use the effective criterion to construct a new class of ballistic diffusions. We show (see Theorem 3.1) that when $d \geq 4$, for any $\eta > 0$, Brownian motion perturbed with a random drift $b(\cdot, \omega)$ such that

$$(1.16) \quad \sup_{x \in \mathbb{R}^d, \omega \in \Omega} |b(x, \omega)| \leq \varepsilon \quad \text{and} \quad \mathbb{E}[b(0, \omega) \cdot e_1] \geq \varepsilon^{2-\eta} \quad \text{for } \varepsilon > 0,$$

satisfies the effective criterion with $\ell = e_1$ if $\varepsilon$ is small enough. The conditions (1.16) allow for laws $\mathbb{P}$ of the environment such that (1.13) does not apply. Indeed, since the constant $c_e$ is larger than 1, as one can see from an inspection of the proof of [11], Theorem 2.5, (1.13) requires that $\mathbb{E}[b(0, \omega) \cdot e_1]$ is larger than $(c_e - 1)\mathbb{E}[(b(0, \omega) \cdot e_1)_-]$ which can be chosen to be of order $\varepsilon$ under (1.16). Note that in the discrete setting, Sznitman (see [20]) obtained similar results under conditions significantly weaker than (1.16). Indeed, he showed that a discrete version of the effective criterion is satisfied by randomly perturbed simple random walk on $\mathbb{Z}^d$ with a drift $d(0, \omega) \stackrel{\text{def}}{=} E_{0,\omega}[X_1 - X_0]$ of size $\varepsilon$ such that $\mathbb{E}[d(0, \omega) \cdot e_1]$ is larger than $\varepsilon^{5/2-\eta}$ when $d = 3$, respectively larger than $\varepsilon^{3-\eta}$ when $d \geq 4$. The strength of this result in contrast to ours is that it includes expected drifts of an order not larger than $\varepsilon^2$, which enabled him to construct examples for condition $(T')$ where Kalikow's condition (see, e.g., [20], (5.3) therein) fails. Considering condition (5.23) of [11] as a continuous analogue of Kalikow's condition, we believe that such examples also exist in our setting. Since, however, the continuous setup with the finite range dependence tends to complicate the arguments, we did not attempt to retrieve the full strength of Sznitman's result.

Let us now briefly describe the proof leading to the new example. In order to verify the effective criterion (1.14) under (1.16), we slice a large box $B$ [as defined below (1.14)] into thinner slabs transversal to $e_1$ and propagate good controls on the exit behavior out of these slabs to the box $B$ using a refinement of the estimate (see Lemma 2.3 and Proposition 3.3) that was instrumental in the induction argument leading to the effective criterion. The heart of the matter is then to prove these good controls for the thinner slabs. To this end, we express the probability that the trajectory exits through the right side of a slab with the help of the Green operator of the diffusion killed when exiting the slab; see (3.23). This quantity is linked to the Green operator of killed Brownian motion via a certain perturbation equality; see (3.40). For Brownian motion, however, an explicit formula obtained by the well-known "method of the images" from electrostatics [see (3.30)] allows us to compute all necessary estimates.

Let us finally explain how this article is organized. In Section 2, we first introduce some notation and then we show the equivalence between the effective criterion and condition $(T')$ when $d \geq 2$; see Theorem 2.6. The key estimate for the induction step is given by Proposition 2.2. In Section 2.1, we discuss the one-dimensional case. In Section 3, we use the effective criterion to show that a certain perturbed Brownian motion satisfies condition $(T')$ when $d \geq 4$. In Section 3.1, we state the main Theorem 3.1 and a refinement of Lemma 2.3; see Proposition 3.3. In Section 3.2, we define the Green operators and Green's functions for which we provide certain deterministic estimates in the case of Brownian motion; see Lemmas 3.7 and 3.9. We also prove a perturbation equality; see Proposition 3.8. In Section 3.3, we use the results from the previous sections to prove the main Theorem 3.1. In Appendix A.1, we give the proof of Lemma 2.3 which is similar to that of [19], Proposition 1.2. In Appendix A.2, we prove Lemma 3.9 using a technique similar to [20], Lemma 2.1.

*Convention on constants.* Unless otherwise stated, constants only depend on the quantities $d, \bar{K}, \nu, R$. We denote with $c$ positive constants with values changing from place to place and with $c_0, c_1, \ldots$ positive constants with values fixed at their first appearance. Dependence on additional parameters appears in the notation.

**2. An effective criterion when $d \geq 2$.** In this section we show that condition $(T')$ [see (1.11)] is equivalent to the effective criterion [see (2.53)] which controls the exit probability from some finite box. By an induction argument we propagate this control to larger boxes that tend to look like infinite slabs. Then one can infer the fast decay of exit probabilities from slabs through "the left" side as required by condition $(T')$.

We first need some notation. For $A, B \subset \mathbb{R}^d$ an open and a closed set, we denote with $T_A = \inf\{t \geq 0; X_t \notin A\}$ the exit time from $A$ and with $H_B = \inf\{t \geq 0; X_t \in B\}$ the entrance time into $B$. For any stopping time $S$, we call $S_0 = 0$, $S_{k+1} = S \circ \theta_{S_k} + S_k$, $k \geq 0$, the iterates of $S$. Here, $\theta_t$ denotes the canonical time shift. We consider a direction $\ell \in S^{d-1}$ and a rotation $\mathcal{R}$ of $\mathbb{R}^d$ such that $R(e_1) = \ell$. The vectors $e_i$, $i = 1, \ldots, d$, constitute the canonical basis. As a shortcut notation for the stopping times in (1.9), we write $T_u = T_u^\ell$ and $\tilde{T}_u = \tilde{T}_u^\ell$, $u \in \mathbb{R}$. Moreover, we introduce

(2.1) $$|z|_\perp = \max_{j \geq 2} |z \cdot \mathcal{R}(e_j)| \qquad \text{for } x \in \mathbb{R}^d.$$

For positive numbers $L, L', \tilde{L}$, we introduce the box

(2.2) $$B = B(\mathcal{R}, L, L', \tilde{L}) \stackrel{\text{def}}{=} \mathcal{R}((-L, L') \times (-\tilde{L}, \tilde{L})^{d-1}),$$

and the positive, respectively negative, part of its boundary

(2.3) $$\partial_+ B = \partial B \cap \{x \in \mathbb{R}^d : \ell \cdot x = L'\}, \qquad \partial_- B = \partial B \setminus \partial_+ B.$$

We also define the following random variables: for $\omega \in \Omega$,

(2.4) $$p_B(\omega) = P_{0,\omega}[X_{T_B} \in \partial_+ B] = 1 - q_B(\omega),$$

(2.5) $$\rho_B(\omega) = \frac{q_B(\omega)}{p_B(\omega)} \in [0, \infty].$$

In the sequel we will use different length scales $L_k, \tilde{L}_k \geq 0$, $k = 0, 1, \ldots$, and the following shortcut notation [cf. (1.5) for the definition of $R$]:

(2.6) $$\begin{aligned} B_k &= B(\mathcal{R}, L_k - R - 1, L_k + 1, \tilde{L}_k) \qquad \text{for } k \geq 0, \\ p_k &= p_{B_k}, \qquad q_k = q_{B_k}, \qquad \rho_k = \rho_{B_k}. \end{aligned}$$

Finally let us set for $k \geq 0$,

(2.7) $$N_k = \frac{L_{k+1}}{L_k}, \qquad n_k = \lfloor N_k \rfloor, \qquad \tilde{N}_k = \frac{\tilde{L}_{k+1}}{\tilde{L}_k}.$$

We start with an easy lemma, introducing the counterpart of a discrete ellipticity constant.

LEMMA 2.1. *Let $C_L$ be the tube $\{z \in \mathbb{R}^d : -\frac{1}{4} < z \cdot e_1 < L, \sup_{2 \leq j \leq d} |z \cdot e_j| < \frac{L}{4}\}$. There exists a constant $0 < \kappa \leq \frac{1}{2}$, such that for any $L \geq 1$, $\omega \in \Omega$, and any rotation $\mathcal{R}$,*

(2.8) $$\begin{aligned} P_{0,\omega}[T_L^{\mathcal{R}(e_1)} < T_{\mathcal{R}(C_L)}] &\geq \kappa^{L+1} \quad \text{and} \\ P_{0,\omega}[\tilde{T}_{-L}^{\mathcal{R}(e_1)} < T_{\mathcal{R}(-C_L)}] &\geq \kappa^{L+1}. \end{aligned}$$

PROOF. We define the function $\psi(s) = \frac{5}{4}\mathcal{R}(e_1)s$, for $0 \leq s \leq 1$. With the support theorem (see [2], page 25), we obtain that there is a constant $c > 0$ such that for all $x \in \mathbb{R}^d, \omega \in \Omega$, $P_{x,\omega}[\sup_{0 \leq s \leq 1} |X_s - X_0 - \psi(s)| < \frac{1}{4}]$ and $P_{x,\omega}[\sup_{0\leq s\leq 1}|X_s - X_0 + \psi(s)| < \frac{1}{4}]$ are both larger than $c$. Then we set $\kappa = \min\{c, \frac{1}{2}\}$. The claim follows by applying the Markov property $\lceil L \rceil$ times. □

We are now ready to prove the main induction step which in essence bounds moments of $\rho_1$ in terms of moments of $\rho_0$.

PROPOSITION 2.2. *($d \geq 2$) There exist $c_1 > R+2, c_2, c_3 > 1$, such that when $N_0 \geq 3$, $L_0 \geq c_1$, $\tilde{N}_0 \geq 150 N_0$, $\tilde{L}_0 \geq R+2$, for any $a \in (0,1]$:*

$$\mathbb{E}[\rho_1^a] \leq c_2 \Bigg\{ \kappa^{-10L_1} \bigg( c_3 \tilde{L}_1^{(d-2)} \frac{L_1^3}{L_0^2} \tilde{L}_0 \mathbb{E}[q_0] \bigg)^{\tilde{N}_0/(12N_0)}$$
$$(2.9)$$
$$+ \sum_{0 \leq m \leq n_0 + 1} (c_3 \tilde{L}_1^{(d-1)} \mathbb{E}[\rho_0^{2a}])^{(n_0+m-1)/2} \Bigg\}.$$

PROOF. For $i \in \mathbb{Z}$ and $L_0 > R+2$, we introduce the slabs of width $R$:

$$(2.10) \qquad \mathcal{S}_i = \left\{ x \in \mathbb{R}^d : iL_0 - \frac{R}{2} \leq x \cdot \ell \leq iL_0 + \frac{R}{2} \right\}$$

and denote by $I(\cdot)$ the function on $\mathbb{R}^d$ such that $I(x) = i$ if $x \cdot \ell - iL_0 \in [-\frac{L_0}{2}, \frac{L_0}{2})$, $i \in \mathbb{Z}$. In particular, $I$ takes the value $i$ on $\mathcal{S}_i$, for all $i \in \mathbb{Z}$. We define the successive times of visit to the different slabs $\mathcal{S}_i$ as the iterates $V_k, k \geq 0$, of the stopping time

$$(2.11) \qquad V = \inf\{t \geq 0 : X_t \in \mathcal{S}_{I(X_0)-1} \cup \mathcal{S}_{I(X_0)+1}\}.$$

We also need the stopping time

$$(2.12) \qquad \tilde{T} = \inf\{t \geq 0 : |X_t|_\perp \geq \tilde{L}_1\}.$$

In a first step we obtain a control on $\mathbb{E}[\rho_1^a]$ using the following quantities: for $\omega \in \Omega, i \in \mathbb{Z}$,

$$(2.13) \qquad \hat{\rho}(i,\omega) = \sup\left\{ \frac{\hat{q}(x,\omega)}{\hat{p}(x,\omega)} : x \in \mathcal{S}_i, |x|_\perp < \tilde{L}_1 \right\},$$

where

$$(2.14) \qquad \hat{q}(x,\omega) = P_{x,\omega}[X_{V_1} \in \mathcal{S}_{I(x)-1}] = 1 - \hat{p}(x,\omega).$$

The first step then comes with the following lemma.

LEMMA 2.3. *Under the assumptions of Proposition 2.2,*
$$(2.15) \quad \mathbb{E}[\rho_1^a] \leq \kappa^{-a(L_1+1)} \mathbb{P}[\mathcal{G}^c] + 2 \sum_{0 \leq m \leq n_0+1} \prod_{-n_0+1 < i \leq m} \mathbb{E}[\hat{\rho}(i,\omega)^{2a}]^{1/2},$$

*where*

$$(2.16) \qquad \mathcal{G} = \{\omega \in \Omega : P_{0,\omega}[\tilde{T} \leq \tilde{T}_{-L_1+R+1} \wedge T_{L_1+1}] \leq \kappa^{9L_1}\}.$$

The proof of this lemma is similar to the proof of [19], (2.39) in Proposition 2.1, or [20], Lemma 1.2. For the reader's convenience, we include the argument in Appendix A.1.

We now complete the proof of Proposition 2.2. Except for a few modifications due to the continuous setup, we follow the steps in the proof of [19], Proposition 2.1. We first bound $\mathbb{P}[\mathcal{G}^c]$ in terms of $\mathbb{E}[q_0]$. In the next section, we infer a different bound on this probability; see (3.13). By Chebyshev's inequality, we find that

$$\mathbb{P}[\mathcal{G}^c] \leq \kappa^{-9L_1} P_0[\tilde{T} \leq \tilde{T}_{-L_1+R+1} \wedge T_{L_1+1}], \quad (2.17)$$

and our task is to derive an upper bound on the right-hand side. We introduce for $u \in \mathbb{R}$, $j \geq 2$, the stopping times

$$\sigma_u^{\pm,j} = \inf\{t \geq 0 : \pm X_t \cdot \mathcal{R}(e_j) \geq u\}, \quad (2.18)$$

$$\bar{L} = 2(n_0+2)(\tilde{L}_0+1) + R, \qquad J = \left\lfloor \frac{\tilde{L}_1}{\bar{L}} \right\rfloor. \quad (2.19)$$

Since $\tilde{L}_1 = \tilde{N}_0 \tilde{L}_0 \geq 150 n_0 \tilde{L}_0$, $n_0 \geq 3$ and $\tilde{L}_0 \geq 2 + R$, it follows that $J \geq 15$. On the event $\{\tilde{T} \leq \tilde{T}_{-L_1+R+1} \wedge T_{L_1+1}\}$, $P_0$-a.s., at least one of the projections $|X_t \cdot \mathcal{R}(e_j)|$, $j \geq 2$, reaches the value $J\bar{L}$ before $X_t$ exits the box $B_1$. Hence:

$$P_0[\tilde{T} \leq \tilde{T}_{-L_1+R+1} \wedge T_{L_1+1}] \leq \sum_{j \geq 2} P_0[\sigma_{J\bar{L}}^{+,j} \leq T_{B_1}] + P_0[\sigma_{J\bar{L}}^{-,j} \leq T_{B_1}]. \quad (2.20)$$

Let us write $\sigma_u$ in place of $\sigma_u^{+,2}$ and bound the term $P_0[\sigma_{J\bar{L}} \leq T_{B_1}]$, the other terms being treated similarly. The strong Markov property yields that

$$P_0[\sigma_{J\bar{L}} \leq T_{B_1}] \leq \mathbb{E} E_{0,\omega}[\sigma_{(J-1)\bar{L}} \leq T_{B_1}, P_{X_{\sigma_{(J-1)\bar{L}}},\omega}[\sigma_{J\bar{L}} \leq T_{B_1}]]. \quad (2.21)$$

We define the auxiliary box

$$B' = B(\mathcal{R}, L_0 - R, L_0, \tilde{L}_0 + 1); \quad (2.22)$$

see (2.2) for the notation and let $H^i$, $i \geq 0$, denote the iterates of the stopping time $H^1 = T_{B_1} \wedge T_{X_0 + B'}$. Then for any $\omega \in \Omega$, $x \in B_1$ with $x \cdot \mathcal{R}(e_2) = (J-1)\bar{L}$, we have

$$P_{x,\omega}[\sigma_{J\bar{L}} > T_{B_1}] \geq P_{x,\omega}\left[\bigcap_{k=0}^{2(n_0+1)-1} \theta_{H^k}^{-1}\{H^1 < T_{\partial_- B' + X_0}\}\right], \quad (2.23)$$

because on the event in the right-hand side, the trajectory either exits $B_1$ before $\sigma_{J\bar{L}}$ right away on $\{H^1 < T_{\partial_- B' + X_0}\}$ or it exits the box $B_1$ through "the right," since for every $k \geq 0$, on $\theta_{H^k}^{-1}\{H^1 < T_{\partial_- B' + X_0}\}$ the trajectory $P_{x,\omega}$-a.s. moves between time $H^k$ and $H^{k+1}$ at most a distance $\tilde{L}_0 + 1$ into direction $\mathcal{R}(e_2)$ and at least a distance $L_0$ into direction $\ell$ until it leaves $B_1$, and since

$$(2.24) \qquad 2(n_0+1)(\tilde{L}_0+1) = \bar{L} - 2(\tilde{L}_0+1) - R < \bar{L}$$

and $2(n_0+1)L_0 > 2L_1 - R$, the width of $B_1$ in direction $\ell$. In order to obtain a lower bound on the right-hand side of (2.23) with the help of the strong Markov property, we cover the set $G(J-1) \stackrel{\text{def}}{=} \{x \in B_1 : |x \cdot \mathcal{R}(e_2) - (J-1)\bar{L}| \leq 2(n_0+1)(\tilde{L}_0+1)\}$, which contains the trajectories up to $T_{B_1}$ described by the event in the right-hand side of (2.23), with a collection of disjoint and rotated unit cubes $C_m$ with centers $x_m$. The cardinality of this collection is proportional to the volume of $G(J-1)$.

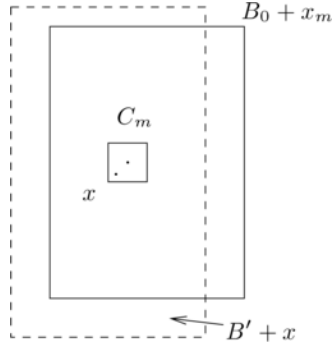

Fig. 1. *Graphical explanation of (2.25).*

For any $k,m \geq 0$ and any $\omega \in \Omega$, we have that on $\{X_{H^k} \in C_m\}$, $P_{0,\omega}$-a.s.,

$$(2.25) \quad P_{X_{H^k},\omega}[X_{T_{B'+X_0}} \in \partial_+ B' + X_0] \geq P_{X_{H^k},\omega}[X_{T_{B_0+x_m}} \in \partial_+ B_0 + x_m],$$

as for any $x \in C_m$, it follows from the definitions of $B'$ [see (2.22)] and $B_0$ [see (2.6)] that $\partial_+ B_0 + x_m \subset (\overline{B'}+x)^c$, $\partial_- B' + x \subset \overline{B_0}^c + x_m$ and $x \in B_0 + x_m$; see Figure 1. Here $\overline{U}$ denotes the closure of $U \subset \mathbb{R}^d$. Therefore any piece of trajectory contained in $B_0 + x_m$, connecting $x \in C_m$ to $\partial_+ B_0 + x_m$ has to exit $B' + x$, but cannot touch $\partial_- B' + x$.

As a consequence, we deduce from (2.23) using the strong Markov property that for any $\omega \in \Omega, x \in B_1$ with $x \cdot \mathcal{R}(e_2) = (J-1)L$,

$$(2.26) \quad P_{x,\omega}[\sigma_{J\bar{L}} > T_{B_1}] \geq \left( \inf_m \inf_{x \in C_m} P_{x,\omega}[X_{T_{B_0+x_m}} \in \partial_+ B_0 + x_m] \right)^{2(n_0+1)}$$
$$\stackrel{\text{def}}{=} 1 - \phi(J-1,\omega),$$

and thus, in view of (2.21), we find

$$(2.27) \qquad P_0[\sigma_{J\bar{L}} \leq T_{B_1}] \leq \mathbb{E}[P_{0,\omega}[\sigma_{(J-2)\bar{L}} \leq T_{B_1}]\phi(J-1,\omega)].$$

From (2.24), we see that $G(J-1) \subset \{x \in B_1 : x \cdot \mathcal{R}(e_2) \geq (J-2)\bar{L} + 2(\tilde{L}_0 + 1) + R\}$, and therefore the random variable $\phi(J-1, \cdot)$ is $\mathcal{H}_{\{z \cdot \mathcal{R}(e_2) \geq (J-2)\bar{L}+R\}}$-measurable whereas $P_{0,\cdot}[\sigma_{(J-2)\bar{L}} \leq T_{B_1}]$ is $\mathcal{H}_{\{z \cdot \mathcal{R}(e_2) \leq (J-2)\bar{L}\}}$-measurable. Thus the finite range dependence property implies that

$$(2.28) \qquad P_0[\sigma_{J\bar{L}} \leq T_{B_1}] \leq P_0[\sigma_{(J-2)\bar{L}} \leq T_{B_1}] \mathbb{E}[\phi(J-1, \omega)].$$

Using the notation (2.26) and observing that $1 - p^k \leq k(1-p)$ for $k \geq 1$, $p \geq 0$, we obtain

$$(2.29) \quad \mathbb{E}[\phi(J-1, \omega)] \leq 2(n_0 + 1)\mathbb{E}\left[\sup_m \sup_{x \in C_m} P_{x,\omega}[X_{T_{B_0+x_m}} \in \partial_- B_0 + x_m]\right].$$

We now observe that the cardinality of the collection of cubes $C_m$ is proportional to $2L_1 \cdot 4(n_0+1)(\tilde{L}_0+1) \cdot (2\tilde{L}_1)^{d-2} \leq c\tilde{L}_1^{d-2} \frac{L_1^2}{L_0} \tilde{L}_0$. Then translation invariance and an application of Harnack's inequality to the harmonic function $x \mapsto P_{x,\omega}[X_{T_{B_0}} \in \partial_- B_0]$ yield that

$$(2.30) \qquad \mathbb{E}[\phi(J-1, \omega)] \leq c'\tilde{L}_1^{d-2} \frac{L_1^3}{L_0^2} \tilde{L}_0 \mathbb{E}[q_0],$$

where we used the notation (2.4). Coming back to (2.28), we see that

$$(2.31) \qquad \begin{aligned} P_0[\sigma_{J\bar{L}} \leq T_{B_1}] &\leq P_0[\sigma_{(J-2)\bar{L}} \leq T_{B_1}] \\ &\quad \times c\tilde{L}_1^{d-2} \frac{L_1^3}{L_0^2} \tilde{L}_0 \mathbb{E}[q_0], \quad \text{and by induction} \\ &\leq \left\{c\tilde{L}_1^{d-2} \frac{L_1^3}{L_0^2} \tilde{L}_0 \mathbb{E}[q_0]\right\}^m \qquad \text{for all } 0 \leq m \leq \left\lfloor \frac{J}{2} \right\rfloor. \end{aligned}$$

Similar bounds hold for each term in the right-hand side of (2.20) and since[1] $\lfloor \frac{J}{2} \rfloor \geq \frac{\tilde{N}_0}{12N_0}$ from our assumptions on $N_0$, $\tilde{N}_0$ and $\tilde{L}_0$ we conclude from (2.17), (2.20) and (2.31) that

$$(2.32) \qquad \mathbb{P}[\mathcal{G}^c] \leq \kappa^{-9L_1} 2(d-1) \left\{c\tilde{L}_1^{d-2} \frac{L_1^3}{L_0^2} \tilde{L}_0 \mathbb{E}[q_0]\right\}^{\tilde{N}_0/(12N_0)}.$$

So far we found an upper bound for the first term of the right-hand side of (2.15). To complete the proof of (2.9), we are now going to bound the second term.

For any $i \in \mathbb{Z}$, we cover the set $\{x \in \mathcal{S}_i : |x|_\perp < \tilde{L}_1\}$ [appearing in the definition of $\hat{\rho}(i, \omega)$; see (2.13)] with a collection of disjoint and rotated unit

---

[1] $\lfloor \frac{J}{2} \rfloor \geq \frac{J}{2} - \frac{1}{2} \geq \frac{\tilde{N}_0 \tilde{L}_0}{4(n_0+2)(\tilde{L}_0+1)+2R} - 1 \geq \frac{\tilde{N}_0 \tilde{L}_0}{4(5/3n_0)(3/2)\tilde{L}_0 + n_0\tilde{L}_0} - 1 \geq \frac{\tilde{N}_0}{11N_0} - 1 \geq \frac{\tilde{N}_0}{N_0}(\frac{1}{11} - \frac{1}{150}) \geq \frac{\tilde{N}_0}{12N_0}$.

cubes $\tilde{C}_k$ with cardinality at most $(R+1)(2\tilde{L}_1+1)^{d-1}$. As a result, for $0 < a < 1$,

$$\mathbb{E}[\hat{\rho}(i,\omega)^{2a}] \leq \sum_k \mathbb{E}\left[\frac{\sup_{x\in\tilde{C}_k}\hat{q}(x,\omega)^{2a}}{\inf_{x\in\tilde{C}_k}\hat{p}(x,\omega)^{2a}}\right]. \tag{2.33}$$

By Harnack's inequality, there is a constant $c \geq 1$ such that

$$\frac{\sup_{x\in\tilde{C}_k}\hat{q}(x,\omega)}{\inf_{x\in\tilde{C}_k}\hat{p}(x,\omega)} \leq c^2\frac{\hat{q}(x_k,\omega)}{\hat{p}(x_k,\omega)} \qquad \text{for every } 1 \leq k, \omega \in \Omega.$$

Moreover, observe that $\hat{q}(x_k,\omega) \leq q_0 \circ t_{x_k}\omega$; see (2.6) for the notation. Using translation invariance, we see that the second term on the right-hand side of (2.15) is less than or equal to

$$2 \sum_{0 \leq m \leq n_0+1} ((R+1)(2\tilde{L}_1+1)^{d-1}c^{4a}\mathbb{E}[\rho_0^{2a}])^{(m+n_0-1)/2}. \tag{2.34}$$

Choosing $c_3 \geq (R+1)3^{d-1}c^4$ sufficiently large completes the proof of Proposition 2.2. $\square$

Similarly to [19], we are going to iterate (2.9) along an increasing sequence of boxes $B_k$, which tend to look like infinite slabs transversal to the direction $\ell$. For the definition of these boxes, we consider

$$u_0 \in (0,1], \qquad v = 8, \qquad \alpha = 240, \tag{2.35}$$

and choose two sequences $L_k, \tilde{L}_k, k \geq 0$, such that

$$\begin{aligned}L_0 \geq c_1, \qquad R+2 \leq \tilde{L}_0 \leq L_0^3, \quad &\text{and for } k \geq 0, \\ L_{k+1} = N_k L_k, \qquad \text{with } N_k = \frac{\alpha}{u_0}v^k \quad &\text{and} \quad \tilde{L}_{k+1} = N_k^3 \tilde{L}_k.\end{aligned} \tag{2.36}$$

As a consequence we see that for $k \geq 0$:

$$L_k = \left(\frac{\alpha}{u_0}\right)^k v^{k(k-1)/2} L_0, \tag{2.37}$$

$$\tilde{L}_k = \left(\frac{L_k}{L_0}\right)^3 \tilde{L}_0. \tag{2.38}$$

LEMMA 2.4. *There exists $c_4 \geq c_1$, such that when for some $L_0 \geq c_4$, $R+2 \leq \tilde{L}_0 \leq L_0^3$, $a_0 \in (0,1]$, $u_0 \in [\kappa^{L_0/d}, 1]$,*

$$\varphi_0 \stackrel{\text{def}}{=} c_3 \tilde{L}_1^{(d-1)} L_0 \mathbb{E}[\rho_0^{a_0}] \leq \kappa^{u_0 L_0}, \tag{2.39}$$

*then for all $k \geq 0$,*

$$\varphi_k \stackrel{\text{def}}{=} c_3 \tilde{L}_{k+1}^{(d-1)} L_k \mathbb{E}[\rho_k^{a_k}] \leq \kappa^{u_k L_k} \qquad \text{with } a_k = a_0 2^{-k}, u_k = u_0 v^{-k}. \tag{2.40}$$

AN EFFECTIVE CRITERION AND A NEW EXAMPLE 13As the proof is purely algebraic and hence identical to the proof of [19], Lemma 2.2, we omit it here. We now use the induction result to control the exit behavior from a slab.

PROPOSITION 2.5. *There exists $c_5 \geq c_4, c_6 > 1$, such that when for some $L_0 \geq c_5$, $R + 2 \leq \tilde{L}_0 \leq L_0^3$,*

$$(2.41) \quad c_6 \left(\log \frac{1}{\kappa}\right)^{3(d-1)} \tilde{L}_0^{(d-1)} L_0^{3(d-1)+1} \inf_{a \in (0,1]} \mathbb{E}[\rho_0^a] < 1,$$

*with $B_0, \rho_0$ as in (2.6), then for some $c > 0$,*

$$(2.42) \quad \limsup_{L \to \infty} L^{-1} \exp\{c(\log L)^{1/2}\} \log P_0[\tilde{T}_{-bL} < T_L] < 0 \quad \text{for all } b > 0.$$

PROOF. In view of (2.37), (2.38), we see that (2.39) is equivalent to

$$(2.43) \quad u_0^{-3(d-1)} \kappa^{-u_0 L_0} c_3 \alpha^{3(d-1)} \tilde{L}_0^{(d-1)} L_0 \mathbb{E}[\rho^{a_0}] \leq 1.$$

The minimum of the function $[\kappa^{L_0/d}, 1] \ni u_0 \mapsto u_0^{-3(d-1)} \kappa^{-u_0 L_0}$ is $c'(L_0 \times \log \frac{1}{\kappa})^{3(d-1)}$, provided $L_0 \geq c_5$. Hence choosing $c_6 = 2c'c_3\alpha^{3(d-1)}$, we can make sure that whenever (2.41) holds, for some $L_0 \geq c_5, R + 2 \leq \tilde{L}_0 \leq L_0^3$, then (2.39) holds for some $a_0 \in (0,1], u_0 \in [\kappa^{L_0/d}, 1]$. By Lemma 2.4, (2.40) holds for all $k \geq 0$. For any $b > 0$, we are now looking for a bound on $P_0[\tilde{T}_{-bL} < T_L]$ when $L$ is large. For every large enough $L$, we can find a unique $k$ with

$$(2.44) \quad L_k < bL \leq L_{k+1}.$$

We then introduce the auxiliary box $B'_k = B(\mathcal{R}, L_k - R, L_k, \tilde{L}_k + 1)$, and use an argument similar to (2.23)–(2.26) to find a lower bound for $P_{0,\omega}[\tilde{T}_{-bL} > T_L]$; that is, we require in essence that the trajectory successively exits certain translates of the box $B'_k$ through the "right" side, $\lfloor \frac{L}{L_k} \rfloor + 1$ times. We therefore cover the set $G' \stackrel{\text{def}}{=} \{x \in \mathbb{R}^d, |x|_\perp \leq (\frac{L}{L_k} + 1)(\tilde{L}_k + 1), x \cdot \ell \in (-bL, L)\}$, playing the role of former set $G(J-1)$, with disjoint and rotated unit cubes $C'_j$ with centers $x'_j$. The cardinality of this collection is at most $m_k \stackrel{\text{def}}{=} ((b+1)L + 1)(2(\frac{L}{L_k} + 1)(\tilde{L}_k + 1) + 1)^{(d-1)}$. For $L$ large, we introduce the event

$$(2.45) \quad \Gamma \stackrel{\text{def}}{=} \left\{\omega \in \Omega : \sup_j \sup_{x \in C'_j} P_{x,\omega}[X_{T_{B_k + x'_j}} \in \partial_- B_k + x'_j] \geq \kappa^{(1/2)u_k L_k}\right\}.$$

Then, for any $\omega \in \Gamma^c$, we obtain by arguments as before that

$$(2.46) \quad P_{0,\omega}[\tilde{T}_{-bL} > T_L] \geq (1 - \kappa^{(1/2)u_k L_k})^{(L/L_k + 1)}.$$

On the other hand, using translation invariance, Harnack's inequality and Chebyshev's inequality, we find that there is a $c > 0$ such that

$$\mathbb{P}[\Gamma] \leq c m_k \kappa^{(-1/2)u_k L_k} \mathbb{E}[q_k]. \qquad (2.47)$$

Since $q_k \leq \rho_k^{a_k}$ and because of (2.40), we obtain that

$$\mathbb{P}[\Gamma] \leq c' \frac{m_k}{c_3 \tilde{L}_{k+1}^{(d-1)} L_k} \kappa^{(1/2)u_k L_k}, \qquad (2.48)$$

and a simple computation using (2.44) and (2.36) shows that for large $L$,

$$\begin{aligned}
\frac{m_k}{\tilde{L}_{k+1}^{(d-1)} L_k} &\leq \frac{((b+1)L+1)(2(L/L_k+1)(\tilde{L}_k+1)+1)^{d-1}}{N_k^{3(d-1)} \tilde{L}_k^{d-1} L_k} \\
&\leq \frac{(1/b+2)L_{k+1} 6^{(d-1)}(L_{k+1}/(bL_k)+1)^{d-1} \tilde{L}_k^{d-1}}{N_k^{3(d-1)} \tilde{L}_k^{d-1} L_k} \\
&\leq c(b)\left(\frac{1}{b} + \frac{1}{N_k}\right)^{d-1} N_k^d N_k^{-3(d-1)} \leq c'(b),
\end{aligned} \qquad (2.49)$$

since $d \geq 2$. As a consequence, we obtain from (2.48) that

$$\mathbb{P}[\Gamma] \leq c(b) \kappa^{(1/2)u_k L_k}. \qquad (2.50)$$

Assembling (2.46), (2.50) and using $1 - p^m \leq m(1-p)$, $p, m \geq 0$, we see that for large $L$:

$$P_0[\tilde{T}_{-bL} < T_L] \leq \left(c(b) + \frac{L}{L_k} + 1\right) \kappa^{1/2 u_k L_k}. \qquad (2.51)$$

From (2.40), (2.36), we obtain $u_k L_k = \frac{u_0^2}{\alpha} v^{-2k} L_{k+1} \overset{(2.44)}{\geq} \frac{u_0^2}{\alpha} v^{-2k} bL$, and so the right-hand side of (2.51) is less than $c'(b) N_k \kappa^{1/2(u_0^2/\alpha)v^{-2k}bL}$. From the inequality $L_k \leq bL$ and (2.37), we deduce that if $L$ is large, then $k \leq c(\log bL)^{1/2}$, and we obtain

$$P_0[\tilde{T}_{-bL} < T_L] \leq c'(b) \exp\left\{-\frac{1}{4} \frac{u_0^2}{\alpha} \left(\log \frac{1}{k}\right) bL \exp(-c(\log(bL))^{1/2})\right\}, \qquad (2.52)$$

for large $L$. This implies the claim (2.42). □

We are now ready to prove the main result of this section.

THEOREM 2.6. *There exists a constant $c_7(d) > 1$, such that for $\ell \in S^{d-1}$ the following conditions are equivalent:*

(i) *There exist $a \in (0,1]$ and a box $B = B(\mathcal{R}, L - R - 2, L + 2, \tilde{L})$ with $\mathcal{R}(e_1) = \ell, L \geq c_5, R + 2 \leq \tilde{L} < L^3$ with*

$$(2.53) \qquad c_7 \left(\log \frac{1}{\kappa}\right)^{3(d-1)} \tilde{L}^{d-1} L^{3(d-1)+1} \mathbb{E}[\rho_B^a] < 1,$$

(ii) $(T')$ *holds with respect to $\ell$ [see (1.11)],*
(iii) $(T)_\gamma$ *holds with respect to $\ell$ for some $\gamma \in (\frac{1}{2}, 1)$ [see (1.10)].*

PROOF. The implication (i) implies (ii) is proved in the same way as the corresponding statement in [19], Theorem 2.4. Indeed, we define $c_7 = 2^{(d-1)} c_6$ and observe that as a result of (2.53),

$$(2.54) \qquad c_6 \left(\log \frac{1}{\kappa}\right)^{3(d-1)} \tilde{L}'^{(d-1)} L^{3(d-1)+1} \mathbb{E}[\rho_B^a] < 1,$$

with $\tilde{L}' = (\tilde{L} + 2) \wedge L^3 \in (\tilde{L}, 2\tilde{L})$. If $B'$ denotes the box $B(\mathcal{R}', L - R - 1, L + 1, \tilde{L}')$ and if the rotation $\mathcal{R}'$ is close enough to $\mathcal{R}$,

$$(2.55) \qquad p_B \leq p_{B'} \quad \text{and hence} \quad \rho_{B'} \leq \rho_B.$$

As a result, whenever $\mathcal{R}'$ is sufficiently close to $\mathcal{R}$, we can apply Proposition 2.5 to the box $B'$, and find that

$$(2.56) \qquad \limsup L^{-\gamma} \log P_0[\tilde{T}^{\ell'}_{-bL} < T^{\ell'}_L] < 0$$

for any $\gamma \in (0,1), b > 0$ with $\ell' = \mathcal{R}(e_1)$.

This proves (ii). It is plain that (iii) follows from (ii).

We now show that (iii) implies (2.53). The neighborhood appearing in the definition of $(T)_\gamma$ contains for some small $\alpha > 0$ and all $j = 2, \ldots, d$, the vectors $\ell'_j = \cos(\alpha)\ell + \sin(\alpha)\mathcal{R}(e_j)$, $\ell''_j = \cos(\alpha)\ell - \sin(\alpha)\mathcal{R}(e_j)$. For large $L'$ and $0 < b < 1$, we choose $L + 2 = L' \frac{1-b}{2\cos(\alpha)}$ and $\tilde{L} = L' \frac{1+b}{2\sin(\alpha)}$. (In particular $\tilde{L} \leq L^3$ if $L'$ is large enough depending on $\alpha$ and $b$.) As a consequence, if we set $B = B(\mathcal{R}, L - 2 - R, L + 2, \tilde{L})$, then $\partial_- B$ is included in the region where $x \cdot \ell'_j \leq -bL'$ or $x \cdot \ell''_j \leq -bL'$ for some $2 \leq j \leq d$ (see also Figure 2). In other words,

$$(2.57) \qquad \mathbb{E}[q_B] \leq P_0\left[\text{there exists } \ell' \in \bigcup_{j=2}^d \{\ell'_j, \ell''_j\} : \tilde{T}^{\ell'}_{-bL'} < T^{\ell'}_{L'}\right],$$

and from (1.10), we see that for some $c > 0$:

$$(2.58) \qquad \mathbb{E}[q_B] \leq 2(d-1) e^{-cL^\gamma} \qquad \text{if $L$ is large enough.}$$

Hence for large $L$, for $a \in (0,1)$ and $c' > 0$:

$$(2.59) \qquad \mathbb{E}[\rho_B^a] \leq \mathbb{E}[\rho_B^a, p_B \geq e^{-c'L^\gamma}] + \mathbb{E}[\rho_B^a, p_B < e^{-c'L^\gamma}],$$

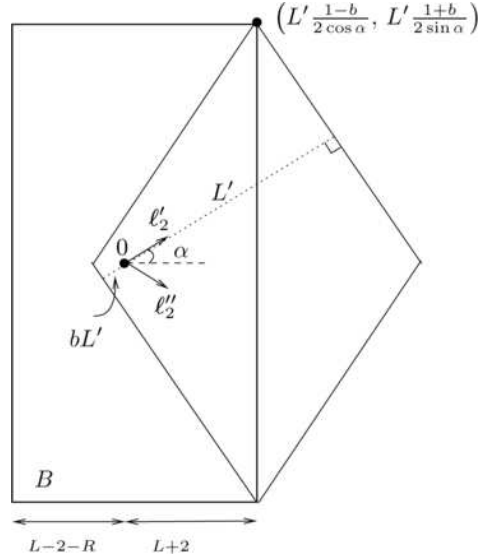

Fig. 2.

so that using the definition (2.5) and Jensen's inequality to bound the first term and $\rho_B \leq \kappa^{-(L+3)}$ because of Lemma 2.1 to control the second term, we find for large $L, a \in (0,1)$ and $c' > 0$:

$$
\begin{aligned}
(2.60) \quad \mathbb{E}[\rho_B^a] &\leq e^{c'aL^\gamma}\mathbb{E}[q_B]^a + \kappa^{-a(L+3)}\mathbb{P}[q_B \geq 1 - e^{-c'L^\gamma}] \\
&\stackrel{(2.58)}{\leq} (2(d-1))^a e^{(c'-c)aL^\gamma} + 2\kappa^{-a(L+3)}2(d-1)e^{-cL^\gamma}.
\end{aligned}
$$

If we choose $a = L^{-1/2}$ and $c' > 0$ sufficiently small, we obtain

$$(2.61) \quad \limsup_{L\to\infty} L^{-(\gamma-1/2)} \log \mathbb{E}[\rho_B^{L^{-1/2}}] < 0.$$

This implies (2.53) and thus finishes the proof of Theorem 2.6. □

REMARK 2.1. As mentioned in the Introduction, it is conjectured that the conditions $(T)$, $(T')$ and $(T)_\gamma$ for a $\gamma \in (0,1)$ are all equivalent and Theorem 2.6 proves part of it. An improvement of the rather crude bound on the second term on the right-hand side of (2.59) is likely to yield the equivalence of $(T')$ and $(T_\gamma)$ also for $\gamma$ smaller than $1/2$. Moreover, the latter theorem together with Proposition 2.5 strengthen the belief that $(T)$ and $(T')$ are equivalent. Indeed, we have in fact obtained the equivalence of $(T')|\ell$ and

$$(*) \quad \limsup_{L\to\infty} L^{-1}\exp\{c(\log L)^{1/2}\}\log P_0[\tilde{T}^{\ell'}_{-bL} < T^{\ell'}_L] < 0 \qquad \text{for all } b > 0$$

and $\ell'$ close to $\ell$ [cf. (2.42)], which is just slightly weaker than $(T)|\ell$, since $\exp\{c(\log L)^{1/2}\}$ grows more slowly than any polynomial. Also note that using [11], (3.36) therein, (∗) actually holds for $\ell' \in S^{d-1}$ satisfying $\ell' \cdot v > 0$, where $v = \lim_{t\to\infty} \frac{X_t}{t}$ denotes the limiting velocity which has been shown to be deterministic and nonzero (ballistic behavior) under $(T')|\ell$ when $d \geq 2$; see [11, 13] and also [19] in the discrete setting.

2.1. *The one-dimensional case.* We introduce here the one-dimensional counterpart of the effective criterion and show that condition $(T)$ is equivalent to $(T')$ and to $P_0$-a.s. transience; see Proposition 2.7. Unlike the multidimensional case, condition $(T')$ does not imply ballistic behavior when $d = 1$, since one can construct one-dimensional diffusions in random environments that tend to infinity, hence satisfy $(T')$, and have zero limiting velocity. A natural question to ask is then whether directional transience, that is, convergence to $\infty$ into some direction, or at least ballistic behavior implies $(T')$ also in higher dimensions.

We first adapt the definitions (2.4), (2.5), (2.10), (2.13) to the one-dimensional setting. Instead of boxes or slabs, we now consider intervals. For any $L > 0$, $\rho_B$ [see (2.5)] is replaced by

$$\rho_L = \frac{P_{0,\omega}[\tilde{T}_{-L} < T_L]}{P_{0,\omega}[\tilde{T}_{-L} > T_L]}. \tag{2.62}$$

For $L_0 \geq 1, i \in \mathbb{Z}$, we redefine $\mathcal{S}_i$ [see (2.10)] as $\mathcal{S}_i = iL_0$. The definition of the stopping times $V_k, k \geq 0$ [see (2.11)], remains unchanged. Then we set for $\omega \in \Omega$, $i \in \mathbb{Z}$:

$$\hat{\rho}(i,\omega) = \frac{\hat{q}(iL_0,\omega)}{\hat{p}(iL_0,\omega)},$$

where

$$\hat{q}(x,\omega) = P_{x,\omega}[X_{V_1} \in \mathcal{S}_{I(x)-1}] = P_{x,\omega}[\tilde{T}_{-L_0+x} < T_{L_0+x}] = 1 - \hat{p}(x,\omega).$$

PROPOSITION 2.7. *(d = 1) The following conditions are equivalent:*

(i) *There exists $a \in (0,1]$, $L > R$, such that $\mathbb{E}[\rho_L^a] < 1$.*
(ii) *There exists $L > R$, such that $\mathbb{E}[\log \rho_L] < 0$.*
(iii) *Condition $(T)$ holds relative to $e_1$.*
(iv) *Condition $(T')$ holds relative to $e_1$.*
(v) $\lim_{t\to\infty} X_t \cdot e_1 = \infty$, $P_0$-*a.s.*

PROOF. The fact that (i) implies (ii) follows from Jensen's inequality since by Lemma 2.1 $\mathbb{E}[\rho_L^a] \leq \kappa^{-a(L+1)} < \infty$. We now show that (ii) implies (iii). We have from (ii) that $-\mu \stackrel{\text{def}}{=} \mathbb{E}[\log \rho_{L_0}] < 0$ for some $L_0 > R$. We are going to use a similar argument as in Appendix A.1 or as in [18], Proposition 2.6 therein. For any $b > 0$ and any real $L > 4L_0/b$, we define

$$(2.63) \qquad n' = \left\lfloor \frac{bL}{L_0} \right\rfloor \quad \text{and set} \quad n_0 = \left\lfloor \frac{L}{L_0} \right\rfloor.$$

[In the spirit of Appendix A.1, $-n'$ plays the role of $-n_0 + 1$; see, e.g., (A.2) or (A.3).] We define the function $f$ on $\{-n', -n'+1, \ldots, n_0+2\}$ by (A.1) and modify the definition of $\tau$ [see (A.3)] as follows:

$$(2.64) \qquad \tau = \inf\{k \geq 0; X_{V_k} \in \mathcal{S}_{n_0+2} \cup \mathcal{S}_{-n'}\}.$$

Since $-\hat{p}(X_{V_m}) + \hat{q}(X_{V_m})\rho(I(X_{V_m}))^{-1}$ vanishes $P_{0,\omega}$-a.s. for all $m \geq 0$, we obtain by an argument similar to (A.4)–(A.7), that for all $\omega \in \Omega$ and $L > \frac{4L_0}{b}$:

$$(2.65) \quad \begin{aligned} & P_{0,\omega}[\tilde{T}_{-bL} < T_L] \\ & \leq \frac{f(0)}{f(-n')} \frac{\prod_{-n',n_0+1}^{-1}}{\prod_{-n',n_0+1}^{-1}} = \frac{\prod_{-n',0}^{-1} + \prod_{-n',1}^{-1} + \cdots + \prod_{-n',n_0+1}^{-1}}{1 + \prod_{-n',-n'+1}^{-1} + \cdots + \prod_{-n',n_0+1}^{-1}} \leq 1. \end{aligned}$$

We then take the expectation with respect to $\mathbb{P}$ of the left-hand side and split it according to the sets where $\sup_{0 \leq k \leq n_0+1} \prod_{-n',k}^{-1}$ is smaller, respectively larger, than $\frac{1}{n_0+2} e^{-c_\mu L}$ with $c_\mu \stackrel{\text{def}}{=} \frac{\mu b}{8L_0}$. As a consequence:

$$(2.66) \quad \begin{aligned} & P_0[\tilde{T}_{-bL} < T_L] \\ & \leq e^{-c_\mu L} + \left(\frac{L}{L_0} + 2\right) \\ & \qquad \times \sup_{0 \leq k \leq n_0+1} \mathbb{P}\left[\sum_{j=-n'+1}^{k} \log \hat{\rho}(j) \geq -c_\mu L - \log(n_0+2)\right]. \end{aligned}$$

Then we decompose the sum appearing in the second term into three sums of independent random variables $\hat{\rho}(j), j = i \mod 3$ where $i = 0, 1$ or $2$. Moreover, since $n_0 \leq \frac{n'+1}{b}$ and by the choice of $c_\mu$, we observe for all $n'$ large enough that for any $0 \leq k \leq n_0 + 1$, we have $\frac{1}{n'+k}(c_\mu L + \log(n_0+2)) \leq \mu/4$. Hence the probability on the right-hand side of (2.66) is less than

$$(2.67) \qquad \sum_{i=0}^{3} \mathbb{P}\left[\frac{3}{n'+k} \sum_{\substack{-n'+1 \leq j \leq k \\ j = i \mod 3}} (-\log \hat{\rho}(j) - \mu) \leq -\mu/4\right].$$

As for any $j \in \mathbb{Z}, \omega \in \Omega$, $|\log \hat{\rho}(j)| \leq (L_0 + 1)\log(\frac{1}{\kappa})$, by (2.8), it follows from an Azuma-type inequality (see, e.g., [1]) that for any $0 \leq k \leq n_0 + 1$, (2.67)

is less than

$$\leq \sum_{i=0}^{3} \exp\left\{-\frac{1}{2}\left(\frac{\mu}{4}\frac{n'+k}{3}\right)^2 |\{j \in [-n'+1,k] \mid j = i \bmod 3\}|^{-1}\right.$$

(2.68)
$$\left. \times ((L_0+1)\log \kappa^{-1} + \mu)^{-2}\right\}$$

$$\leq 3\exp\left\{-c(\mu,L_0)\left(\frac{bL}{L_0} - 1\right)\right\}.$$

In view of (2.66), this implies condition $(T)$; see (1.12).

The implication (iii) $\Rightarrow$ (iv) is clear. To show (iv) $\Rightarrow$ (i), we follow the argument of the corresponding multidimensional statement [see Theorem 2.6, (iii) $\Rightarrow$ (i)], that is, in place of (2.58) and (2.60), we have $P_0[\tilde{T}_{-L} < T_L] \leq e^{-cL^\gamma}$ and $\mathbb{E}[\rho_L^a] \leq e^{(c'-c)aL^\gamma} + 2\kappa^{-a(L+3)}e^{-cL^\gamma}$, for $L$ large.

We now come to the implication (v) $\Rightarrow$ (ii). We follow the arguments in [4], Theorem 2, point (b) in the case of a line. This theorem applied to the discrete Markov chain $X_{V_k}$, $k \geq 0$, under $P_{0,\omega}$ for an $L_0 > R$ in fact shows the equivalence of (ii) and (v). For the reader's convenience, we extract and present here the ideas which are relevant for the implication (v) $\Rightarrow$ (ii). For $\omega \in \Omega$, $L_0 > R$, $n \in \mathbb{Z}$, we introduce the shortcut notation $p_n = \hat{p}(nL_0,\omega) = 1 - q_n$ and $\delta_n \stackrel{\text{def}}{=} P_{nL_0,\omega}[\tilde{T}_{(n-1)L_0} = \infty] \stackrel{\text{def}}{=} 1 - \eta_n$. We claim that

(2.69) $$\mathbb{P}[\delta_0 > 0] = 1.$$

Indeed, let us assume by contradiction that there exists some $\omega'$ in the set of full measure $\{\omega \in \Omega : P_{-L_0,\omega}[X_t \to \infty] = 1\}$ such that $\delta_0(\omega') = 0$. Then repeated use of the strong Markov property shows that $P_{0,\omega'}[\liminf_t X_t \leq -L_0] = 1$, a contradiction.

Next, we see that for $n \in \mathbb{Z}$,

(2.70)
$$\eta_n = P_{nL_0,\omega}[\tilde{T}_{(n-1)L_0} < \infty]$$
$$= q_n + p_n \eta_{n+1} \eta_n \quad \text{and thus} \quad \eta_n = \frac{q_n}{1 - p_n \eta_{n+1}}.$$

As a consequence, for all $\omega \in \Omega, n \leq -1$,

(2.71)
$$\delta_n = 1 - \eta_n \stackrel{(2.70)}{=} 1 - \frac{1-p_n}{1-p_n\eta_{n+1}} = \hat{\rho}(n)^{-1}\eta_n \delta_{n+1}, \quad \text{and by induction}$$
$$= (\hat{\rho}(n)\hat{\rho}(n+1)\cdots\hat{\rho}(-1))^{-1}\eta_n\eta_{n+1}\cdots\eta_{-1}\delta_0.$$

Taking the logarithm of the latter expression and splitting the resulting sum into sums of i.i.d. random variables [similarly as below (2.66)], we obtain from the law of large numbers:

(2.72) $$\lim_{n \to \infty} \frac{1}{n} \log \delta_n = \mathbb{E}[-\log \hat{\rho}(0)] + \mathbb{E}[\log \eta_0], \qquad \mathbb{P}\text{-a.s.},$$

since $\lim_n \frac{1}{n} \log \delta_0 = 0$, $\mathbb{P}$-a.s. by (2.69).

On the other hand, by translation invariance of $\mathbb{P}$, we see that for any $\varepsilon > 0$,

$$\text{(2.73)} \quad \mathbb{P}[|\log \delta_n| > \varepsilon n] = \mathbb{P}[|\log \delta_0| > \varepsilon n] \overset{n \to \infty}{\longrightarrow} \mathbb{P}[\delta_0 = 0] \overset{(2.69)}{=} 0.$$

In other words, $\frac{1}{n} \log \delta_n$ converges to 0 in probability, so the right-hand side of (2.72) vanishes and $\mathbb{E}[\log \hat{\rho}(0)] = \mathbb{E}[\log \eta_0]$ which is strictly negative because of (2.69). This proves the implication (v) $\Rightarrow$ (ii).

To show the converse implication, we use the fact that (ii) implies condition $(T)$. Following [11] [see the proof of $(3.1) \Rightarrow (3.2)$ therein], we observe that $P_0[T_L = \infty] \leq P_0[\tilde{T}_{-L} < T_L]$, since $P_0[\tilde{T}_{-L} = T_L = \infty] = 0$ as in every time unit, the trajectory can escape from the interval $[-L, L]$ with a probability bounded away from 0. Observe that the left-hand side increases with $L$ while the right-hand side tends to 0 by condition $(T)|e_1$. Hence $P_0$-a.s., $\limsup_{t \to \infty} X_t = \infty$. From the strong Markov property and translation invariance of $\mathbb{P}$, we obtain for any $L > 0$:

$$\text{(2.74)} \quad P_0[\tilde{T}_{L/2} \circ \theta_{T_L} < T_{4L/3} \circ \theta_{T_L}] = P_0[\tilde{T}_{-L/2} < T_{5L/6}].$$

Under condition $(T)|e_1$, the right-hand side decreases exponentially and hence an application of Borel–Cantelli's lemma yields that $P_0$-a.s. for large integer $L$, $T_{4L/3} < \tilde{T}_{L/2} \circ \theta_{T_L} + T_L$. As a result, we can $P_0$-a.s. construct an integer-valued sequence $L_k \uparrow \infty$, with $L_{k+1} = \lfloor \frac{4}{3} L_k \rfloor$ and $T_{L_{k+1}} < \tilde{T}_{L_k/2} \circ \theta_{T_{L_k}} + T_{L_k}$, $k \geq 0$. This shows (v). $\square$

REMARK 2.2. Let us mention that for any $L > 0$,

$$\text{(2.75)} \quad \mathbb{E}[\log \rho_L] = -2L \mathbb{E}[b(0)/a(0)],$$

and as a consequence, if conditions (i) or (ii) above are satisfied for some $L \geq R$, they are in fact satisfied for all $L > 0$. Indeed, using the scale function $s(x, \omega) = \int_0^x \exp\{-\int_0^y 2b(u, \omega)/a(u, \omega) \, du\} \, dy$, for $x \in \mathbb{R}$, $\omega \in \Omega$ (see, e.g., [2], pages 78 and 88), we can write $\rho_L = \frac{s(L)}{-s(-L)}$. It follows that for $L > 0$, $\mathbb{E}[\log \rho_L]$ equals

$$\text{(2.76)} \quad \begin{aligned} &\mathbb{E}\left[\log \int_0^L e^{-\int_0^y 2b(u,\omega)/a(u,\omega) \, du} \, dy\right] \\ &\quad - \mathbb{E}\left[\log \int_{-L}^0 e^{-\int_{-L}^y 2b(u,\omega)/a(u,\omega) \, du} e^{\int_{-L}^0 2b(u,\omega)/a(u,\omega) \, du} \, dy\right]. \end{aligned}$$

Because of translation invariance of $\mathbb{P}$, the second term becomes

$$\mathbb{E}\left[\log \int_0^L e^{-\int_0^y 2b(u,\omega)/a(u,\omega) \, du} \, dy\right] + \mathbb{E}\left[\int_{-L}^0 \frac{2b(u,\omega)}{a(u,\omega)} \, du\right],$$

so that the first term of (2.76) is canceled out. Fubini's theorem then yields $\mathbb{E}[\log \rho_L] = -\int_{-L}^{0} \mathbb{E}[\frac{2b(u,\omega)}{a(u,\omega)}]\,du$, and the claim follows from translation invariance of $\mathbb{P}$.

## 3. An example of a ballistic diffusion.

3.1. *Main result and preliminaries.* In this section, we use the effective criterion to show that a Brownian motion perturbed by a small random drift which is bounded by $\varepsilon > 0$ and whose expectation in direction $\ell = e_1$ is of order $\varepsilon^{2-\eta}$ with $\eta > 0$ satisfies condition $(T')|e_1$. The interest of this class of diffusions stems from the fact that it contains new examples of ballistic diffusions which in particular do not fulfill the criterion of [11], Theorem 5.2 therein, which states that there exists a constant $c_e > 1$ such that if

$$(3.1) \qquad \mathbb{E}[(b(0,\omega) \cdot e_1)_+] > c_e \mathbb{E}[(b(0,\omega) \cdot e_1)_-],$$

then $(T)|e_1$ holds. Before we give further explanations on this matter (see Remark 3.1 below), we introduce the family of perturbed Brownian motions studied in this section. For any $\varepsilon \in (0, \bar{K}]$, $\eta > 0$ and $\omega \in \Omega$, we consider the class of diffusions attached to an operator of the form

$$(3.2) \qquad \mathcal{L} = \tfrac{1}{2}\Delta + b(x,\omega) \cdot \nabla,$$

where we require that for all $x \in \mathbb{R}$, $\omega \in \Omega$,

$$(3.3) \qquad |b(x,\omega)| \leq \varepsilon, \qquad \lambda \stackrel{\mathrm{def}}{=} \mathbb{E}[b(0,\omega) \cdot e_1] \geq \varepsilon^{2-\eta}.$$

Note that the constant $\bar{K}$, the ellipticity constant $\nu$ and the dependence range $R$ [see (1.2)–(1.5)] do not depend on $\varepsilon$. We keep the convention concerning constants stated at the end of the Introduction. Moreover, when we write that an expression holds "for large enough $L$" we mean that the expression holds for all $L$ larger than some $c(\eta)$.

The main result of the section is

THEOREM 3.1. *When $d \geq 4$, for any $\eta \in (0,1)$ there is $\varepsilon_0(\eta, d) > 0$ such that whenever (3.3) holds for $0 < \varepsilon \leq \varepsilon_0$, then condition $(T')|e_1$ is satisfied.*

REMARK 3.1. Clearly, (3.1) is equivalent to

$$(3.4) \qquad \mathbb{E}[b(0,\omega) \cdot e_1] > (c_e - 1)\mathbb{E}[(b(0,\omega) \cdot e_1)_-].$$

An inspection of the proof of Theorem 5.2 in [11], reveals that $c_e > 1$, and hence (3.4) fails when $\varepsilon > 0$ is small, if $\mathbb{E}[b(0,\omega) \cdot e_1]$ is of order $\varepsilon^{2-\eta}$ with $0 < \eta < 1$ and $\sup_{\omega \in \Omega}(b(0,\omega) \cdot e_1)_-$ is of order $\varepsilon$ under an adequate choice of $\mathbb{P}$. With this observation one can rather straightforwardly produce examples where (3.3) holds with $\varepsilon < \varepsilon_0(\eta, d)$, but (3.1) or (3.4) fails.

The rest of the section is devoted to the proof of Theorem 3.1. We will verify the effective criterion (2.53) when $\varepsilon$ is smaller than some $\varepsilon_0(\eta, d)$ for $a = 1/2$ and a box $B = B(\text{Id}, NL' - R - 2, NL' + 2, \frac{1}{4}(NL')^3)$ [see (2.2)], where

$$(3.5) \quad N = L^3 \quad \text{and} \quad L = L' - \frac{R}{2} \text{ is an integer such that } L = \left\lfloor \frac{1}{4\varepsilon} \right\rfloor.$$

The starting point to estimate $\mathbb{E}[\rho_B^{1/2}]$ is (2.15). Here we set [cf. (2.6)] $L_1 = NL', \tilde{L}_1 = \frac{1}{4}(NL')^3$, $L_0 = L'$, $n_0 = N$ and $a = 1/2$. With these choices, the box $B$ defined above, on which we want to check (2.53), equals $B_1 + e_1$. In order to apply (2.15) we use the following.

LEMMA 3.2. *For $a \in (0, 1]$ and $B_1$ a box as in (2.6) with $\ell = e_1, L_1 \geq R + 3$ and $\mathcal{R} = \text{Id}$,*

$$(3.6) \quad \mathbb{E}[\rho_{B_1+e_1}^a] \leq c^a \mathbb{E}[\rho_{B_1}^a].$$

PROOF. Since for every $\omega \in \Omega, P_{x,\omega}[X_{T_{B_1}} \in \partial_\pm B_1]$ is harmonic on $(-2, 2)^d$, Harnack's inequality implies that

$$\frac{P_{-e_1,\omega}[X_{T_{B_1}} \in \partial_- B_1]}{P_{-e_1,\omega}[X_{T_{B_1}} \in \partial_+ B_1]} \leq c\rho_{B_1}(\omega).$$

The claim then follows from translation invariance of $\mathbb{P}$. □

For the purpose of this section, we need a bound on $\mathbb{P}[\mathcal{G}^c]$ appearing in (2.15) which differs from (2.32) and which is essentially the same as the estimate in [20], Theorem 1.1. We now follow [20] to introduce the notation used for this bound. Let $h, H, M$ be positive integers with

$$(3.7) \quad 2h \leq H \leq \frac{(NL')^3}{32} \quad \text{and} \quad M = \left\lfloor \frac{(NL')^3}{32H} \right\rfloor.$$

Later on [see (3.51)], we will choose $H$ and $h$ to be of order $(NL')^2$ and $L^2$, respectively. We introduce the exit time $S$ from a tube:

$$(3.8) \quad S = \inf\left\{t > 0; |(X_t - X_0) \cdot e_1| \geq L \text{ or } \sup_{j \geq 2} |(X_n - X_0) \cdot e_j| \geq h\right\}$$

and the expected displacement

$$(3.9) \quad \Delta(x, \omega) = E_{x,\omega}[X_S] - x, \qquad x \in \mathbb{R}^d, \omega \in \Omega.$$

Moreover, for $0 < \gamma \leq 1$, later chosen to be of order $\varepsilon^{1-\eta}$ [see (3.51)], we define

$$(3.10) \quad p_L = \inf_{j \geq 2} \mathbb{P}[\text{for all } z \in \tilde{B}^j, \Delta(z, \omega) \cdot e_1 \geq \gamma L],$$

where for $2 \leq j \leq d$,

(3.11) $$\tilde{B}^j = \{y \in B, |y \cdot e_j| < H\}.$$

Let us now state the analogue of Theorem 1.1 in [20].

PROPOSITION 3.3. *There exists a constant $c_8 > R + 3$ such that when $L \geq c_8$ and*

(3.12) $$\delta^{-1} \stackrel{\text{def}}{=} \exp\left\{-\frac{\gamma N}{128}\right\} + \frac{10N}{\gamma} \exp\left\{-\frac{\gamma N}{32}\left(\frac{H}{2hN} - \frac{4}{\gamma}\right)_+^2\right\} < 1,$$

*then for any $0 < a \leq 1$*

(3.13) $$\mathbb{E}[\rho_B^a] \leq c^a \kappa^{-aNL'} 2d \exp\left\{-\frac{M}{2}\left(p_L - \frac{10NL}{M}\frac{\log \kappa^{-1}}{\log \delta}\right)_+^2\right\}$$
$$+ c^a \frac{2\mathbb{E}[\hat{\rho}(0,\omega)^{2a}]^{N/2}}{(1 - \mathbb{E}[\hat{\rho}(0,\omega)^{2a}]^{1/2})_+}.$$

Since the proof is very similar to the one of Theorem 1.1 in [20], we only make a few comments here. Because of (3.6), we can estimate $\mathbb{E}[\rho_B^a]$ with the help of (2.15). We bound the second term on the right-hand side in the latter expression using translation invariance of $\mathbb{P}$ and obtain the second term on the right-hand side of (3.13). The intuitive idea behind the estimate on $\mathbb{P}[\mathcal{G}^c]$ in the first term on the right-hand side of (2.15), leading to the first term on the right-hand side of (3.13), is to consider nested boxes $\hat{B}_k = (-NL' + R + 2, NL' + 2) \times (-k4H, k4H)^{d-1}$ for $0 \leq k \leq M$ contained in the big box $B$. Then in order not to exit through "the left or right" of $B$, the trajectory has to reach the boundary of box $\hat{B}_k$ before exiting $B$ and then move from box $\hat{B}_k$ to box $\hat{B}_{k+1}$ without exiting $B$. The probability of this last step is related to the quantity $1 - p_L$.

Note that the coefficient in the first term of $\delta^{-1}$ differs from the result in [20] as the width of $B$ is a multiple of $L'$ while the definition of the time $S$ uses the quantity $L = L' - \frac{R}{2}$. This affects the right-hand side of the expression below (1.24) in [20].

Despite the finite range dependence, the remark in [20] below (1.29) still holds since (in the notation of [20]) the random variables $Z_k(e)$ and $Z_{k-1}(e)$ are measurable respectively in $\mathcal{H}_{\{z \in \mathbb{R}^d : z \cdot e \geq 4kH - H - h\}}$ and in $\mathcal{H}_{\{z \in \mathbb{R}^d : z \cdot e \leq 4(k-1)H + H + h\}}$. The involved half-spaces are separated by a distance $2(H - h)$ which is larger than $H$ by (3.7). Hence $(Z_k)_{0 \leq k \leq M}$ are independent if $H$ and thus $L$ are large enough.

3.2. *Bounds on the Green operator.* The main Theorem 3.1 will follow after choosing $h$ and $H$ as in (3.51) once we show exponential decay in $L \propto \varepsilon^{-1}$ of both terms on the right-hand side of (3.13) for $a = 1/2$. Therefore the goals of this section are to find a tractable expression for $\hat{\rho}(0, \omega)$ (see Lemma 3.5) that involves the Green operator of the diffusion killed when exiting the open slab $\mathcal{S} \stackrel{\text{def}}{=} \{x \in \mathbb{R}^d : |x \cdot e_1| < L\}$, and then investigate its relation with the Green operator of killed Brownian motion; see Proposition 3.8. Certain deterministic estimates on the latter operator and its kernel (see Lemmas 3.7 and 3.9) will then be instrumental in the proof of the desired exponential decay of $\mathbb{E}[\rho_B^{1/2}]$; see Proposition 3.10.

Throughout this section, we use the shortcut notation $b_1 \stackrel{\text{def}}{=} b \cdot e_1$ and we set $\|f\|_\infty = \sup_{x \in \mathcal{S}} |f(x)|$, for any function $f$ on $\mathcal{S}$. For any bounded measurable function $f$ on $\mathcal{S}$ and any $x \in \mathcal{S}$, $\omega \in \Omega$, let us denote with

$$
\begin{aligned}
G_\mathcal{S}^\omega f(x) &\stackrel{\text{def}}{=} E_{x,\omega}\left[\int_0^{T_\mathcal{S}} f(X_s)\, ds\right], \qquad \text{respectively} \\
G_\mathcal{S} f(x) &\stackrel{\text{def}}{=} E\left[\int_0^{T_\mathcal{S}} f(x + W_s)\, ds\right],
\end{aligned}
\tag{3.14}
$$

the Green operator of the diffusion, respectively Brownian motion, killed when exiting the slab $\mathcal{S}$. (Here $E$ denotes the expectation with respect to some measure under which $W_s$ is a Brownian motion.) Note that by (3.16) below, these operators acting on $L^\infty$ have norm bounded by $2L^2$. Moreover, the semi-group $P_t^\omega$ of the diffusion in environment $\omega$ killed when exiting $\mathcal{S}$ is defined as

$$
(3.15) \qquad P_t^\omega f(x) = E_{x,\omega}[f(X_t), t < T_\mathcal{S}] \qquad \text{for } x \in \mathcal{S}, \quad t \geq 0.
$$

In a similar fashion, we denote with $P_t$ the semi-group of a Brownian motion killed when exiting $\mathcal{S}$.

The following lemma states basic bounds on the expected exit time from the slab $\mathcal{S}$ and on the supremum-norm of the operator $P_t^\omega$.

LEMMA 3.4. *For $\omega \in \Omega, \varepsilon \in (0, 1/4)$, $x \in \mathcal{S}$, under the assumption (3.3) and with the definition (3.5),*

$$
(3.16) \qquad \tfrac{2}{3}(L^2 - (x \cdot e_1)^2) \leq E_{x,\omega}[T_\mathcal{S}] \leq 2(L^2 - (x \cdot e_1)^2).
$$

*For any bounded measurable function $f$ and any $\omega \in \Omega$,*

$$
(3.17) \qquad \|P_t^\omega f\|_\infty \leq c_{10} \|f\|_\infty \exp(-c_{11} t / L^2) \qquad \text{for } t > 0.
$$

PROOF. To show (3.16), we consider for $x \in \mathcal{S}, \omega \in \Omega$ the $P_{x,\omega}$-martingale

$$
(3.18) \quad (X_{t \wedge T_\mathcal{S}} \cdot e_1)^2 - (X_0 \cdot e_1)^2 - \int_0^{t \wedge T_\mathcal{S}} 2b_1(X_s, \omega)(X_s \cdot e_1)\, ds - t \wedge T_\mathcal{S}.
$$

After taking expectations and using the monotone convergence theorem, we obtain (3.16) from our assumption $|b(\cdot,\cdot)| \leq \varepsilon$ [see (3.3)] and the choice $L \leq \frac{1}{4\varepsilon}$ [see (3.5)].

We now turn to (3.17). By the support theorem (see [2]) applied to the rescaled diffusion $\frac{1}{L}X_{L^2 t}$ and the fact that $|Lb_1| \leq \frac{1}{4}$, the probability under $P_{x,\omega}$ that the trajectories leave the slab within time $L^2$ when starting in $x \in \mathcal{S}$ is bounded away from 0 by some constant $c_{11}$. Hence the strong Markov property yields for any $t > 0$, $x \in \mathcal{S}, \omega \in \Omega$, that $P_{x,\omega}[t \leq T_{\mathcal{S}}] \leq c_{10} \exp(-c_{11} t/L^2)$, and (3.17) follows from the definition (3.15). □

REMARK 3.2.

1. For Brownian motion starting at $x \in \mathcal{S}$, the expected exit time from the slab $\mathcal{S}$ equals $L^2 - (x \cdot e_1)^2$. The analogue of (3.17) for Brownian motion is also valid.
2. We point out that since $T_{\mathcal{S}}$ has a finite moment under $P_{x,\omega}$ by (3.16), Fubini's theorem applied to (3.14) yields for any bounded measurable function $f$ and any $\omega \in \Omega, x \in \mathcal{S}$ that

$$(3.19) \qquad G^\omega_{\mathcal{S}} f(x) = \int_0^\infty P^\omega_t f(x)\, dt.$$

Of course, the same relation holds for the killed Brownian motion.

Let us now introduce the following shortcut notation for the set appearing in the definition of $\hat{\rho}(0,\omega)$ [see (2.13)]:

$$(3.20) \qquad \mathcal{V} \stackrel{\text{def}}{=} \left\{ x \in \mathbb{R}^d; |x \cdot e_1| \leq \frac{R}{2}, |x|_\perp \leq \frac{1}{4}(NL')^3 \right\}.$$

For later purposes, we observe that (3.16) and our assumption (3.3) on $\lambda$ imply that there are constants $c_{12} > 0$ and $L_1(c_{12}, \eta)$ such that when $L \geq L_1$, then for any $x \in \mathcal{V}, \omega \in \Omega$,

$$(3.21) \qquad G^\omega_{\mathcal{S}} \lambda(x) = \lambda E_{x,\omega}[T_{\mathcal{S}}] \geq c_{12} L^\eta.$$

The next lemma provides a tractable expression of $\rho(0,\omega)$ in terms of the Green operator $G^\omega_{\mathcal{S}}$.

LEMMA 3.5. *For $L \geq 3R, \omega \in \Omega$, with (3.3) and (3.5),*

$$(3.22) \qquad \hat{\rho}(0,\omega) = \sup_{x \in \mathcal{V}} \frac{L - x \cdot e_1 - G^\omega_{\mathcal{S}}(b_1(\cdot,\omega))(x)}{L + x \cdot e_1 + G^\omega_{\mathcal{S}}(b_1(\cdot,\omega))(x)} \leq 5.$$

*[See (2.13), (3.20) for the notation.]*

PROOF. For any $x \in \mathcal{S}$, $\omega \in \Omega$, $X_{t \wedge T_\mathcal{S}} \cdot e_1 - X_0 \cdot e_1 - \int_0^{t \wedge T_\mathcal{S}} b_1(X_s, \omega) \, ds$ is a $P_{x,\omega}$-martingale. Hence, after taking expectations, we obtain from the dominated convergence theorem that [see (2.14) for the notation]

$$\hat{p}(x,\omega) = \frac{x \cdot e_1 + L + G_\mathcal{S}^\omega(b_1(\cdot,\omega))(x)}{2L}. \qquad (3.23)$$

Inserting this expression into the definition (2.13) of $\hat{\rho}(0,\omega)$ yields the claimed equality. Using (3.16), (3.5), we see that for all $L > 0$,

$$|G_\mathcal{S}^\omega(b_1(\cdot,\omega))(x)| \leq \frac{L}{2}, \qquad (3.24)$$

and thus the inequality in (3.22) follows when $L \geq 3R$. □

In order to explore the relationship between the Green operators of the diffusion and Brownian motion [see (3.40)], we need to collect a few facts about the semi-group of Brownian motion. From [16], Theorem 8.1.18, we have that whenever $f$ is a continuous and bounded function, then $(t,x) \mapsto P_t f(x)$ is bounded and in $C^{1,2}([0,\infty) \times \mathcal{S}, \mathbb{R})$. Moreover,

$$\frac{\partial}{\partial t} P_t f = \frac{1}{2} \Delta P_t f \qquad \text{in } (0,\infty) \times \mathcal{S}, \qquad (3.25)$$

$$\lim_{t \to 0} P_t f(x) = f(x), \qquad x \in \mathcal{S}. \qquad (3.26)$$

Since every point on the boundary of $\mathcal{S}$ is regular according to [16], (8.1.16) therein, we have the following continuity property at the boundary (see [16], Theorem 8.1.18):

$$\lim_{\substack{(t,x) \to (s,a) \\ (t,x) \in (0,\infty) \times \mathcal{S}}} P_t f(x) = 0 \qquad \text{for } (s,a) \in (0,\infty) \times \partial \mathcal{S}. \qquad (3.27)$$

Our next step is to express $P_t$ and $G_\mathcal{S}$ in terms of kernels using "the method of images" from electrostatics.

PROPOSITION 3.6. *Let $f$ be a bounded measurable function on $\mathbb{R}^d$. If we define for $t > 0; x, y \in \mathcal{S}$*

$$p(t,x,y) = \sum_{k=-\infty}^{\infty} p_d(t,x,y+2k2Le_1) - p_d(t,x,y^* + (2k+1)2Le_1), \qquad (3.28)$$

*where $p_d(t,x,y) \stackrel{\text{def}}{=} (2\pi t)^{-d/2} \exp\{|x-y|^2/2t\}$ is the d-dimensional heat kernel and $y^*$ is the image of $y$ under reflection with respect to $\{z \in \mathbb{R}^d : z \cdot e_1 = 0\}$, then*

$$P_t f(x) = \int_\mathcal{S} p(t,x,y) f(y) \, dy. \qquad (3.29)$$

*Moreover, when $d \geq 4$, if we define Green's function for distinct $x, y \in \mathcal{S}$ by*

$$(3.30) \quad g(x,y) \stackrel{\text{def}}{=} \sum_{k=-\infty}^{\infty} g_d(x, y + 2k2Le_1) - g_d(x, y^* + (2k+1)2Le_1),$$

*where $g_d(x,y) \stackrel{\text{def}}{=} \int_0^\infty p_d(t,x,y)\,dt = \gamma_d |x-y|^{2-d}$ for $x \neq y$ and an appropriate constant $\gamma_d$, then*

$$(3.31) \quad G_\mathcal{S} f(x) = \int_\mathcal{S} g(x,y) f(y)\,dy.$$

PROOF. The fact that $p(t,x,y)$ in (3.28) satisfies the equality in (3.29) follows from [8], Proposition 8.10, after mapping the interval $[0,a]$ to $[-L,L]$ and after multiplying with $p_{d-1}$. It is well known that $g_d(x,y)$ equals $\gamma_d |x-y|^{2-d}$ for an appropriate constant $\gamma_d$ when $d \geq 3$ and $x \neq y$ (see, e.g., [16], (8.4.10)). To see that the expression in (3.30) is indeed the kernel of $G_\mathcal{S}$, we observe that $p(t,x,y)$ is integrable over $t$ for $x \neq y$, since by the monotone convergence theorem, we have $\int_0^\infty p(t,x,y)\,dt \leq \sum_{k=-\infty}^{\infty} g_d(x, y + 2k2Le_1) + g_d(x, y^* + (2k+1)2Le_1)$, and since the latter series converges absolutely when $d \geq 4$. Moreover, with dominated convergence,

$$(3.32) \quad g(x,y) = \int_0^\infty p(t,x,y)\,dt \qquad \text{for } x \neq y.$$

Then we insert (3.29) into (3.19) and since $(t,y) \mapsto p(t,x,y)f(y)$ is product integrable by Tonelli's theorem and (3.17), we obtain (3.31) from Fubini's theorem and (3.32). □

The next lemma provides gradient estimates on the semi-group and the Green operator of killed Brownian motion which play an important role in the derivation of the perturbation equality (3.40) and in the proof of Proposition 3.10.

LEMMA 3.7. $(d \geq 4)$ *For any bounded, continuous function $f$, there exist $c_{13}, c_{14} > 0$ such that for all $x \in \mathcal{S}$, $t > 0$ and $L > 0$,*

$$(3.33) \quad |\nabla P_t f(x)| \leq \left(\frac{c_{13}}{L} + \frac{c_{14}}{\sqrt{t}}\right) \exp\left(-\frac{c_{11}}{2} t/L^2\right) \|f\|_\infty,$$

$$(3.34) \quad |\nabla G_\mathcal{S} f(x)| \leq c_{15} \|f\|_\infty L.$$

PROOF. We first show (3.33). Let $(x^{(1)}, \ldots, x^{(d)})$ denote the coordinates of a point $x$ in $\mathbb{R}^d$. We estimate the partial derivatives $\partial_i$, $i = 1, \ldots, d$, of $P_t f(x)$ separately. As a consequence of the semi-group property, we have that for $t > 0$, $x \in \mathcal{S}$,

$$(3.35) \quad P_t f(x) = \int_\mathcal{S} p(t/2, x, z) P_{t/2} f(z)\,dz.$$

We let $U \subset \mathcal{S}$ be a neighborhood of $x$. To compute $\partial_i P_t(x)$ by interchanging derivation and integration, we need to show that $|\partial_i p(t/2, x, z) P_{t/2} f(z)|$ is $dx \times dz$ integrable over $U \times \mathcal{S}$. After an application of (3.17), we see that

$$\sup_{x \in U} \int_{\mathcal{S}} |\partial_i p(t/2, x, z) P_{t/2} f(z)| \, dz$$
(3.36)
$$\leq \exp(-c_{11} t / L^2) \|f\|_\infty \sup_{x \in U} \int_{\mathcal{S}} |\partial_i p(t/2, x, z)| \, dz.$$

For $i = 1$, according to (3.28), $|\partial_1 p(t/2, x, z)|$ is smaller than

$$\prod_{j=2}^{d} p_1(t/2, x^{(j)}, z^{(j)}) \sum_{k=-\infty}^{\infty} |\partial_1 p_1(t/2, x^{(1)}, z^{(1)} + 4kL)|$$
(3.37)
$$+ |\partial_1 p_1(t/2, x^{(1)}, -z^{(1)} + (2k+1)2L)|.$$

The integral over $\mathbb{R}^d$ of the first $d-1$ factors in the latter expression equals 1 and using monotone convergence, we find that for any $x \in U$, $t > 0$, the integral on the right-hand side of (3.36) is smaller than

$$\sum_{k \neq 0} p_1(t/2, x^{(1)}, L + 4kL) + p_1(t/2, x^{(1)}, -L + 4kL)$$

$$+ \sum_{k \neq -1} p_1(t/2, x^{(1)}, L + (2k+1)2L)$$

(3.38)
$$+ \sum_{k \neq 0} p_1(t/2, x^{(1)}, -L + (2k+1)2L)$$

$$+ \int_{-L}^{L} \frac{1}{\sqrt{\pi t}} e^{(-1/t)(x^{(1)} - z)^2} \left| \frac{x^{(1)} - z}{t} \right| dz$$

$$+ p_1(t/2, x^{(1)}, -L) + p_1(t/2, x^{(1)}, L).$$

For any $x \in U$, the function $z \mapsto p_1(t/2, x^{(1)}, z)$ is monotone on $(-\infty, -2L]$ and on $[2L, \infty)$. Therefore the first sum in (3.38) is less than $\frac{4}{L} \int_{-\infty}^{\infty} p_1(t/2, x^{(1)}, z) \, dz = \frac{4}{L}$. A similar argument yields that the second and third sums in (3.38) are less than $\frac{c}{L}$. The integral in (3.38) is less than $2 \int_0^\infty \frac{1}{\sqrt{\pi t}} e^{-u^2} u \, du = \frac{1}{\sqrt{\pi t}}$ and the last two terms can also be bounded by $\frac{c}{\sqrt{t}}$. Collecting our estimates, we obtain for $i = 1$ that the left-hand side of (3.36) is less than

(3.39)
$$\exp(-c_{11} t / L^2) \|f\|_\infty \left( \frac{c}{L} + \frac{c'}{\sqrt{t}} \right).$$

Hence we can interchange the derivative $\partial_1$ with the integral in (3.35), and for any $x \in \mathcal{S}$, $|\partial_1 P_t(x)|$ is bounded by (3.39). Similar bounds on $|\partial_i P_t(x)|$,

for $2 \leq i \leq d$, follow from an easier version of the above arguments. Indeed, in an expression corresponding to (3.37), the last factor containing the sum, which was more delicate to treat, will not be affected by the derivative $\partial_i$ and thus its integral over $[-L, L]$ equals 1. This proves (3.33). Since the latter estimate shows that $\nabla P_t f(x)$ is integrable with respect to $t > 0$, (3.34) is an immediate consequence of (3.19). $\square$

The link between the Green operators of the killed diffusion and the killed Brownian motion is expressed by the following perturbation equality:

PROPOSITION 3.8. *Let $f$ be a bounded, continuous function on $\mathbb{R}^d$. Then we have for all $x \in \mathcal{S}, \omega \in \Omega$ that*

$$(3.40) \qquad G_{\mathcal{S}}^{\omega} f(x) = G_{\mathcal{S}} f(x) - G_{\mathcal{S}}^{\omega}(b(\cdot, \omega) \cdot \nabla) G_{\mathcal{S}} f(x).$$

PROOF. The classical idea of the proof is to take the derivative of $P_t^{\omega} P_{u-t} f(x)$ with respect to $t$, which yields $P_t^{\omega}((\mathcal{L} - \Delta) P_{u-t} f)(x)$. Then one integrates both sides with respect to $t$ from 0 to $u$ and with respect to $u$ from 0 to infinity. The result then follows from Fubini's theorem. Let us now present the details of the proof. For $\omega \in \Omega, u > 0, x \in \mathcal{S}$, we claim that

$$(3.41) \qquad P_u^{\omega} f(x) - P_u f(x) = \int_0^u P_t^{\omega}(b(\cdot, \omega) \cdot \nabla P_{u-t} f)(x) \, dt.$$

To prove the claim, we define for $h > 0$ the function

$$(3.42) \qquad e(t, x) \stackrel{\text{def}}{=} P_{u+h-t} f(x) \qquad \text{with } 0 \leq t \leq u, x \in \mathcal{S}.$$

According to [16], Theorem 8.1.18, $e$ is in $C^{1,2}((0, u) \times \mathcal{S})$. Hence we can apply Itô's formula to a function $e_n \in C^{1,2}((0, u) \times \mathbb{R}^d)$ such that $e_n(t, \cdot) = e(t, \cdot)$ on $D_n \stackrel{\text{def}}{=} \{x \in \mathcal{S}, \text{dist}(x, \partial \mathcal{S}) \geq 1/n\}$ and $e_n(t, \cdot) = 0$ on $\mathcal{S}^c$. Because of (3.25), we obtain for all $\omega \in \Omega$, $x \in D_n$ after taking expectations:

$$(3.43) \quad E_{x,\omega}\left[e(u \wedge T_{D_n}, X_{u \wedge T_{D_n}}) - e(h \wedge T_{D_n}, X_{h \wedge T_{D_n}}) - \int_{h \wedge T_{D_n}}^{u \wedge T_{D_n}} b(X_s, \omega) \cdot \nabla e(s, X_s) \, ds\right] = 0.$$

When $n$ tends to $\infty$, $t \wedge T_{D_n}$ increases to $t \wedge T_{\mathcal{S}}$, and it follows from the dominated convergence theorem and (3.27) that for any $\omega \in \Omega, x \in \mathcal{S}$,

$$(3.44) \qquad E_{x,\omega}[e(u \wedge T_{D_n}, X_{u \wedge T_{D_n}}), u \geq T_{\mathcal{S}}] \stackrel{n \to \infty}{\longrightarrow} 0.$$

The same result holds for $h$ in place of $u$. From (3.33), we have that $\sup_{0 \leq t \leq u, x \in \mathcal{S}} |\nabla e(v, x)|$ is finite. Thus coming back to (3.43) and letting

$n \to \infty$, we obtain with dominated convergence that for any $x \in \mathcal{S}$,

$$
\begin{aligned}
(3.45) \quad & E_{x,\omega}[e(u, X_u), u < T_\mathcal{S}] - E_{x,\omega}[e(h, X_h), h < T_\mathcal{S}] \\
& = E_{x,\omega}\left[\int_h^{u \wedge T_\mathcal{S}} b(X_s, \omega) \cdot \nabla e(s, X_s)\, ds, h < T_\mathcal{S}\right].
\end{aligned}
$$

We now insert the definition (3.42) into the above expression and let $h$ tend to 0 using dominated convergence. This concludes the proof of (3.41).

The integral with respect to $u > 0$ of the left-hand side of (3.41) equals $G_\mathcal{S}^\omega f(x) - G_\mathcal{S} f(x)$; see (3.19). On the right-hand side, (3.33) and (3.17) imply that the iterated integral

$$
\begin{aligned}
(3.46) \quad & \int_0^\infty \int_0^\infty |P_t^\omega(b(\cdot, \omega) \cdot \nabla P_{u-t} f)|(x) 1_{\{t < u\}}\, du\, dt \\
& \le \int_0^\infty \int_t^\infty c\varepsilon \|f\|_\infty \left(\frac{c_{13}}{L} + \frac{c_{14}}{\sqrt{u-t}}\right) e^{-c_{11}/(2L^2)(t+u)}\, du\, dt
\end{aligned}
$$

is finite. Hence we can integrate the right-hand side of (3.41) with respect to $u$, use Fubini's theorem and then substitute $u - t$ with $u$. It follows for $\omega \in \Omega, x \in \mathcal{S}$,

$$
(3.47) \quad G_\mathcal{S}^\omega f(x) - G_\mathcal{S} f(x) = \int_0^\infty \int_0^\infty P_u^\omega(b(\cdot, \omega) \cdot \nabla P_t f)(x)\, du\, dt.
$$

The same argument as before allows us to interchange the integrals once more. Finally with (3.33) and a further application of Fubini's theorem we can move the $dt$-integral inside $P_u^\omega(\cdot)$ and interchange it with the gradient. This finishes the proof of Proposition 3.8. $\square$

We close this section with estimates on the Green's function (3.30) of killed Brownian motion and on its gradient. They are at the heart of the proof of Proposition 3.10.

LEMMA 3.9. *($d \ge 4$) For all $x, y \in \mathcal{S}$ and $L > 0$ we have*

$$(3.48) \quad g(x, y) \le c_{16}|x-y|^{2-d} \exp(-c_{17}|x-y|_\perp/L),$$

$$(3.49) \quad |\nabla g(x, y)| \le (c_{18}|x-y|^{1-d} + c_{19}L^{1-d}) \exp(-c_{17}|x-y|_\perp/L).$$

*Moreover, for any bounded Hölder continuous function $f$, $G_\mathcal{S} f$ is twice continuously differentiable on $\mathcal{S}$ and*

$$(3.50) \quad \tfrac{1}{2}\Delta G_\mathcal{S} f(x) = -f(x) \qquad \text{for } x \in \mathcal{S}.$$

The proof is included in Appendix A.2 and the arguments showing (3.48) and (3.49) are similar to the proof of [20], (2.11), (2.13) therein.

3.3. *Proof of Theorem* 3.1. The starting point for the proof is (3.13) with $a = \frac{1}{2}$. We first specify the quantities $h, H, \gamma$ involved in the first term on the right-hand side of (3.13) [see (3.7), (3.10)]:

(3.51)
$$h \stackrel{\text{def}}{=} L'^2, \qquad H \stackrel{\text{def}}{=} \lfloor (NL')^2 \rfloor,$$
$$\gamma \stackrel{\text{def}}{=} \tfrac{1}{4} c_{12} L^{\eta-1}.$$

It is clear that the main Theorem 3.1 follows from the effective criterion once we show exponential decay in $L \propto \varepsilon^{-1}$ of both terms on the right-hand side of (3.13). We first examine the second term. It suffices to show that for large enough $L$

(3.52) $$\mathbb{E}[\hat{\rho}(0, \omega)] \leq \exp\left(-\frac{c_{12}}{2} L^{-1}\right),$$

where $c_{12}$ is defined in (3.21). Indeed, since we assumed $N = L^3$ [see (3.5)], the second term of (3.13) then becomes smaller than $cL \exp(-\frac{c_{12}}{4} L^2)$, which will be more than sufficient for the application of the effective criterion (2.53).

To prove (3.52), we use (3.22) and write $\mathbb{E}[\hat{\rho}(0, \omega)]$ as

(3.53)
$$\mathbb{E}\left[\sup_{x \in \mathcal{V}} \frac{L - x \cdot e_1 - G_{\mathcal{S}}^{\omega}(b_1(\cdot, \omega))(x)}{L + x \cdot e_1 + G_{\mathcal{S}}^{\omega}(b_1(\cdot, \omega))(x)}, \inf_{x \in \mathcal{V}} G_{\mathcal{S}}^{\omega}(b_1(\cdot, \omega))(x) \geq \frac{c_{12}}{2} L^{\eta}\right]$$
$$+ 5\mathbb{P}\left[\inf_{x \in \mathcal{V}} G_{\mathcal{S}}^{\omega}(b_1(\cdot, \omega))(x) < \frac{c_{12}}{2} L^{\eta}\right].$$

When $L$ is larger than some $c(\eta)$, the first term becomes smaller than $1 - \frac{c_{12}}{2} L^{\eta-1} \leq \exp(-\frac{c_{12}}{2} L^{\eta-1})$. Hence (3.52) follows from the next proposition which estimates the second term of (3.53).

PROPOSITION 3.10. *($d \geq 4$) For any $\eta \in (0,1)$, under the assumption (3.3) and with (3.5), we have that*

(3.54) $$\limsup_{L \to \infty} L^{-2/3\eta} \log \mathbb{P}\left[\inf_{x \in \mathcal{V}} G_{\mathcal{S}}^{\omega}(b_1(\cdot, \omega))(x) < \frac{c_{12}}{2} L^{\eta}\right] < 0,$$

*where $\mathcal{V}$ and $c_{12}$ are defined in (3.20) and (3.21).*

Before proving the proposition, we show that (3.54) together with our choices in (3.51) also yield exponential decay of the first term on the right-hand side of (3.13), which then finishes the proof of the main theorem. Using (3.5), we find that

(3.55) $$\delta^{-1} \leq \exp(-cL^{2+\eta}) + c'L^{4-\eta} \exp\{-c''L^{2+\eta}(L^3 - c'''L^{1-\eta})^2\},$$

which tends to 0 as $L$ goes to $\infty$, so that (3.12) holds when $L$ is large. If in addition, we know that [see (3.10) for the notation]

$$\liminf_{L\to\infty} p_L = 1, \tag{3.56}$$

an easy calculation using (3.51) and $M \geq c_{13}NL$ [see (3.7) for the definition] shows that for $L$ large enough, the first term on the right-hand side of (3.13) is less than $c\exp(-cNL)$, and the effective criterion (2.53) is satisfied for large $L$.

We now prove that Proposition 3.10 implies (3.56). First we cover the sets $\tilde{B}_j, 2 \leq j \leq d$ [see (3.11)] with a collection of disjoint cubes of side length $\frac{R}{2}$. The cardinality of this collection is for large $L$ at most $L^\nu$ where $\nu$ only depends on $d$. Translation invariance then yields

$$p_L \geq 1 - \sup_{2\leq j\leq d} c'L^\nu \mathbb{P}\bigg[\inf_{x\in[-R/2,R/2]^d} \Delta(x,\omega)\cdot e_1 < \gamma L\bigg]. \tag{3.57}$$

In this expression we will in essence replace $\Delta(x,\omega)\cdot e_1$ with $G_{\mathcal{S}}^\omega(b_1(\cdot,\omega))(x)$. More precisely, we claim that for large $L$ and for all $\omega \in \Omega$, $x \in [-\frac{R}{2},\frac{R}{2}]^d$,

$$|\Delta(x,\omega)\cdot e_1 - G_{\mathcal{S}}^\omega(b_1(\cdot,\omega))(x)| \leq c_{20}. \tag{3.58}$$

Then with our choice of $\gamma$ [see (3.51)] and with (3.57), Proposition 3.10 implies (3.56) since $[-\frac{R}{2},\frac{R}{2}]^d \subset \mathcal{V}$. We now prove (3.58). The martingale argument leading to (3.23) also shows that for any $x \in \mathcal{S}$, $\omega \in \Omega$

$$G_{\mathcal{S}}^\omega(b_1(\cdot,\omega))(x) = E_{x,\omega}[X_{T_\mathcal{S}}\cdot e_1] - x\cdot e_1. \tag{3.59}$$

The support theorem (see [2]) applied to the rescaled diffusion $\frac{1}{L}X_{L^2t}$ yields a lower bound $c > 0$ (uniform in $x \in \mathbb{R}^d$, $\omega \in \Omega$) for the probability under $P_{x,\omega}$ that $X$ exits a cube of side length $L$ centered at $x$ through the "left or right." Hence with the strong Markov property, for all $\omega \in \Omega$, $x \in \mathbb{R}^d$,

$$P_{x,\omega}[S < \tilde{T}_{-L+x\cdot e_1} \wedge T_{L+x\cdot e_1}] \leq 2(d-1)(1-c)^L, \tag{3.60}$$

which becomes smaller than $L^{-1}$ for large enough $L$. Since $|X_S \cdot e_1| \leq L + |x\cdot e_1|$, $P_{x,\omega}$-a.s. we obtain from (3.59) and (3.60) that for large enough $L$ and for all $\omega \in \Omega, x \in [-\frac{R}{2},\frac{R}{2}]^d$, the left-hand side of (3.58) is less than

$$|E_{x,\omega}[(X_S - X_{T_\mathcal{S}})\cdot e_1, S = \tilde{T}_{-L+x\cdot e_1}\wedge T_{L+x\cdot e_1}]| + c. \tag{3.61}$$

On the event $\{S = \tilde{T}_{-L+x\cdot e_1}\wedge T_{L+x\cdot e_1}\}\cap\{(X_S\cdot e_1)(X_{T_\mathcal{S}}\cdot e_1) > 0\}$, the trajectory $P_{x,\omega}$-a.s. leaves the slab $\mathcal{S}$ and the box $[-L,L]\times[-h,h]^{d-1}+x$ "through the same side." Hence on this event, $|(X_S - X_{T_\mathcal{S}})\cdot e_1| \leq \frac{R}{2}, P_{x,\omega}$-a.s. for $x \in [-\frac{R}{2},\frac{R}{2}]^d$. It remains to show that for all $\omega \in \Omega$, $x \in [-\frac{R}{2},\frac{R}{2}]^d$,

$$\begin{aligned}&|E_{x,\omega}[(X_S - X_{T_\mathcal{S}})\cdot e_1,\\&\quad S = \tilde{T}_{-L+x\cdot e_1}\wedge T_{L+x\cdot e_1}, (X_S\cdot e_1)(X_{T_\mathcal{S}}\cdot e_1) < 0]| \leq c.\end{aligned} \tag{3.62}$$

When $x \cdot e_1 = 0$ the above quantity vanishes. We now consider the case where $0 < x \cdot e_1 \leq \frac{R}{2}$. The remaining case is treated analogously. We find that for $0 < x \cdot e_1 \leq \frac{R}{2}$,

$$
\begin{aligned}
(3.63) \quad & P_{x,\omega}[S = \tilde{T}_{-L+x \cdot e_1} \wedge T_{L+x \cdot e_1}, (X_S \cdot e_1)(X_{T_S} \cdot e_1) < 0] \\
& \leq P_{x,\omega}[T_L < \tilde{T}_{-L+x \cdot e_1} < T_{L+x \cdot e_1}] + P_{x,\omega}[\tilde{T}_{-L+x \cdot e_1} < T_L < \tilde{T}_{-L}].
\end{aligned}
$$

We estimate the first term on the right-hand side. The strong Markov property implies that for all $\omega \in \Omega, 0 < x \cdot e_1 \leq \frac{R}{2}$,

$$
\begin{aligned}
(3.64) \quad & P_{x,\omega}[T_L < \tilde{T}_{-L+x \cdot e_1} < T_{L+x \cdot e_1}] \\
& \leq E_{x,\omega}[T_L < \tilde{T}_{-L+x \cdot e_1}, P_{X_{T_L},\omega}[\tilde{T}_{-L+R/2} < T_{L+R/2}]].
\end{aligned}
$$

The function $e(x) \stackrel{\text{def}}{=} -e^{4\varepsilon x \cdot e_1} + e^{4\varepsilon(L+R/2)}$ satisfies $\mathcal{L}e(x) < 0$ since $|b(\cdot,\cdot)| \leq \varepsilon$. Hence $e(X_t)$ is a supermartingale under $P_{x,\omega}$ for any $x \in \mathbb{R}^d, \omega \in \Omega$. Since $e(x)$ is nonnegative when $x \cdot e_1 \leq L + \frac{R}{2}$, Chebyshev's inequality and the stopping theorem yield for any $y \in \mathbb{R}^d$ with $y \cdot e_1 = L$,

$$
(3.65) \quad P_{y,\omega}[\tilde{T}_{-L+R/2} < T_{L+R/2}] \leq \frac{E_{y,\omega}[e(X_{\tilde{T}_{-L+R/2} \wedge T_{L+R/2}})]}{e^{4\varepsilon(L+R/2)} - e^{4\varepsilon(-L+R/2)}}
$$

$$
\leq \frac{1 - e^{-4\varepsilon R/2}}{1 - e^{-8\varepsilon L}} \leq c\varepsilon \leq c'L^{-1},
$$

for large enough $L$. Inserting this bound into (3.64) and repeating the same type of argument for the second term on the right-hand side of (3.63), we obtain that its left-hand side is of order $L^{-1}$. This finishes the proof of (3.62) since $(X_S - X_{T_S}) \cdot e_1$ is of order $L$, $P_{x,\omega}$-a.s. for $x \in [-\frac{R}{2}, \frac{R}{2}]^d$. Thus (3.58) follows in view of (3.61). As a consequence, Proposition 3.10 implies (3.56) and the main theorem follows as we explained below (3.56).

PROOF OF PROPOSITION 3.10. The idea of the proof is to decompose the $e_1$ projection of the drift $b_1(x,\omega)$ into its expectation $\mathbb{E}[b_1 \cdot e_1] = \lambda$ and a mean-zero term $\tilde{b}(x,\omega)$. As a consequence, the Green operator applied to $b_1$ splits into two terms: a leading term $G_S^\omega \lambda$ which is larger than twice the bound imposed on the Green operator in the event of interest in (3.54) by our choice of constants and by (3.21); an error term $G_S^\omega \tilde{b}$ that we decompose using the perturbation equality (3.40) and which turns out to make no substantial contribution to the leading term with high probability. Hence the event of interest in (3.54) is very unlikely. We now give the details of the proof. Let us introduce the box

$$
(3.66) \quad \mathcal{U} \stackrel{\text{def}}{=} \{x \in \mathbb{R}^d; |x \cdot e_1| \leq L - 1, |x|_\perp \leq \tfrac{1}{4}(NL')^3 + L^2\}
$$

which will be useful later in a discretization step where we need to restrict ourselves to points located at a constant distance of $\partial \mathcal{S}$. As mentioned above we define [see (3.3)]

$$\tilde{b} \stackrel{\text{def}}{=} b_1 - \lambda. \tag{3.67}$$

[For the sake of simplicity we drop the $\omega$ dependence of $b_1, \tilde{b}$ from the notation.] Then the perturbation equality (3.40) applied to $G_{\mathcal{S}}^{\omega}\tilde{b}$ together with (3.21) yields that for large enough $L$,

$$
\begin{aligned}
\mathbb{P}\bigg[\inf_{x \in \mathcal{V}} G_{\mathcal{S}}^{\omega} b_1(x) &< \frac{c_{12}}{2} L^{\eta}\bigg] \\
&\leq \mathbb{P}\bigg[\inf_{x \in \mathcal{V}} G_{\mathcal{S}}\tilde{b}(x) - G_{\mathcal{S}}^{\omega}(b \cdot \nabla) G_{\mathcal{S}}\tilde{b}(x) \leq -\frac{c_{12}}{2} L^{\eta}, \\
&\qquad \sup_{y \in \mathcal{U}} |\nabla G_{\mathcal{S}}\tilde{b}(y)| \leq L^{-1+\eta/3}\bigg] \\
&\quad + \mathbb{P}\bigg[\sup_{y \in \mathcal{U}} |\nabla G_{\mathcal{S}}\tilde{b}(y)| > L^{-1+\eta/3}\bigg].
\end{aligned}
\tag{3.68}
$$

The proposition obviously follows once we prove the following three claims: there exist $\nu', \nu'' \geq 1$ depending only on $d$ such that for large enough $L$,

$$
\text{on the set } \bigg\{\omega \in \Omega; \sup_{y \in \mathcal{U}} |\nabla G_{\mathcal{S}}\tilde{b}(y)| \leq L^{-1+\eta/3}\bigg\},
$$
$$\sup_{x \in \mathcal{V}} |G_{\mathcal{S}}^{\omega}(b \cdot \nabla) G_{\mathcal{S}}\tilde{b}(x)| \leq c L^{\eta/3}, \tag{3.69}$$

$$\mathbb{P}\bigg[\sup_{y \in \mathcal{U}} |\nabla G_{\mathcal{S}}\tilde{b}(y)| > L^{-1+\eta/3}\bigg] \leq L^{\nu'} \exp(-c' L^{2/3\eta}), \tag{3.70}$$

$$\mathbb{P}\bigg[\inf_{x \in \mathcal{V}} G_{\mathcal{S}}\tilde{b}(x) \leq -\frac{c_{12}}{4} L^{\eta}\bigg] \leq L^{\nu''} \exp(-c' L^{c_{21}+2\eta}) \tag{3.71}$$

where $c_{21} = 1$ when $d = 4$ and $c_{21} = 2$ when $d \geq 5$.

We now show (3.69). In view of (3.34) and (3.5), we have that $\sup_{x \in \mathcal{S}} |\nabla G_{\mathcal{S}}\tilde{b}(x)| \leq c_{15} 2\varepsilon L \leq c_{15}/2$. Therefore for any $\omega \in \Omega$ satisfying $\sup_{y \in \mathcal{U}} |\nabla G_{\mathcal{S}}\tilde{b}(y)| \leq L^{-1+\eta/3}$ and any $x \in \mathcal{V}$ we find that

$$
\begin{aligned}
|G_{\mathcal{S}}^{\omega}(b \cdot \nabla) G_{\mathcal{S}}\tilde{b}(x)| &\leq \varepsilon L^{-1+\eta/3} G_{\mathcal{S}}^{\omega} 1_{\mathcal{U}}(x) \\
&\quad + \varepsilon \frac{c_{15}}{2} G_{\mathcal{S}}^{\omega} 1_{\{z \in \mathcal{S}; \text{dist}(z, \partial \mathcal{S}) \leq 1\}}(x) \\
&\quad + \varepsilon \frac{c_{15}}{2} G_{\mathcal{S}}^{\omega} 1_{\{z \in \mathcal{S}; |z|_{\perp} \geq 1/4(NL')^3 + L^2\}}(x).
\end{aligned}
\tag{3.72}
$$

The first term on the right-hand side is smaller than $\frac{1}{4}L^{-2+\eta/3}E_{x,\omega}[T_{\mathcal{S}}] \leq \frac{1}{2}L^{\eta/3}$ by (3.16).

To bound the second term on the right-hand side of (3.72), we define for $L \geq 4(1+R)$ the auxiliary set $\hat{\mathcal{S}} = \{x \in \mathcal{S}; \operatorname{dist}(x, \partial \mathcal{S}) < 2\}$. With a martingale argument similar to (3.18), (3.16), we obtain that for any $\omega \in \Omega$ and $x \in \mathcal{S}$, $E_{x,\omega}[T_{\hat{\mathcal{S}}}] \leq (1-2\varepsilon)^{-1} \leq 2$. Then we introduce the successive times of entrance in $\{x \in \mathbb{R}^d; |x \cdot e_1| \geq L-1\}$ and departure from $\{x \in \mathbb{R}^d; |x \cdot e_1| > L-2\}$:

$$R_1 = T_{L-1} \wedge \tilde{T}_{-L+1}, \qquad D_1 = T_{\{x \in \mathbb{R}^d; |x \cdot e_1| > L-2\}} \circ \theta_{R_1} + R_1,$$

(3.73) and by induction for $k \geq 2$,

$$R_k = R_1 \circ \theta_{D_{k-1}} + D_{k-1}, \qquad D_k = T_{\{x \in \mathbb{R}^d; |x \cdot e_1| > L-2\}} \circ \theta_{R_k} + R_k.$$

With the help of these definitions we now express the Green operator appearing in the second term on the right-hand side of (3.72): for any $\omega \in \Omega, x \in \mathcal{V}$, we have

$$G_{\mathcal{S}}^{\omega}(x) 1_{\{z \in \mathcal{S}; \operatorname{dist}(z, \partial \mathcal{S}) \leq 1\}}$$

(3.74)
$$= \sum_{k \geq 1} E_{x,\omega}\left[\int_{R_k}^{D_k \wedge T_{\mathcal{S}}} 1_{\{z \in \mathcal{S}; \operatorname{dist}(z, \partial \mathcal{S}) \leq 1\}}(X_s)\, ds, R_k < T_{\mathcal{S}}\right]$$
$$\leq \sum_{k \geq 1} E_{x,\omega}[E_{X_{R_k},\omega}[T_{\hat{\mathcal{S}}}], R_k < T_{\mathcal{S}}]$$
$$\leq 2 \sum_{k \geq 1} P_{x,\omega}[R_k < T_{\mathcal{S}}].$$

The sum is bounded by a constant since the strong Markov property and the support theorem imply that for $k \geq 1, x \in \mathcal{V}$, $P_{x,\omega}[R_k < T_{\mathcal{S}}] \leq (1-c)^{k-1}$. Hence the second term on the right-hand side of (3.72) is less than $c'L^{-1}$.

We now examine the last term on the right-hand side of (3.72). We call $\tilde{U}$ the set $\{z \in \mathcal{S}; |z|_{\perp} \geq 1/4(NL')^3 + L^2\}$ appearing in that term. For any $\omega \in \Omega, x \in \mathcal{V}$, the Markov property yields

(3.75)
$$G_{\mathcal{S}}^{\omega} 1_{\tilde{U}}(x) = E_{x,\omega}\left[E_{X_{H_{\tilde{U}}},\omega}\left[\int_0^{T_{\mathcal{S}}} 1_{\tilde{U}}(X_s)\, ds\right], H_{\tilde{U}} < T_{\mathcal{S}}\right]$$
$$\leq \sup_{z \in \mathcal{S}} E_{z,\omega}[T_{\mathcal{S}}] P_{x,\omega}[H_{\tilde{U}} < T_{\mathcal{S}}].$$

Using (3.16) and a scaling argument similar to the one leading to (3.60), we find that the latter expression is smaller than $cL^2 e^{-c'L}$. As a consequence, the last term on the right-hand side of (3.72) is smaller than $L^{-1}$ for large enough $L$. This proves (3.69).

Next we turn to the proof of (3.70). In order to deal with the supremum over the set $\mathcal{U}$, we cover $\mathcal{U}$ with disjoint cubes of side-length $\varepsilon^3$ and centers $y_i, i \in \mathcal{I}$, where $|\mathcal{I}| \leq cL^{12d-8}$. If $Q$ is such a cube with center $y_i$, then according to Lemma 3.9, $-\frac{1}{2}G_{\mathcal{S}}\tilde{b}(y)$ is twice continuously differentiable on $Q' \stackrel{\text{def}}{=} y_i + (-\frac{1}{2}, \frac{1}{2})^d \subset \mathcal{S}$ and satisfies the equation $\Delta u = \tilde{b}$ on $Q'$. Therefore [5], (3.20), page 41, applies and we find that for any $y \in Q$

$$
(3.76) \quad \begin{aligned} &|\nabla G_{\mathcal{S}}\tilde{b}(y) - \nabla G_{\mathcal{S}}\tilde{b}(y_i)| \\ &\leq c|y - y_i|\left(\sup_{z \in Q'}|G_{\mathcal{S}}\tilde{b}(z)| + \sup_{z \in Q'}|\tilde{b}(z)|\right)\left(\left|\log\left(\frac{c'}{|y - y_i|}\right)\right| + 1\right). \end{aligned}
$$

Since the bounds in (3.16) also hold for Brownian motion, we have that $\sup_{z \in Q'}|G_{\mathcal{S}}\tilde{b}(z)| \leq 2L^2 2\varepsilon \leq L$. Thus the right-hand side of (3.76) is less than $cL^{-2+\eta}$ for large enough $L$. With this discretization step we obtain for large enough $L$:

$$
(3.77) \quad \begin{aligned} &\mathbb{P}\left[\sup_{y \in \mathcal{U}}|\nabla G_{\mathcal{S}}\tilde{b}(y)| > L^{-1+\eta/3}\right] \\ &\leq \sum_{i \in \mathcal{I}}\mathbb{P}[|\nabla G_{\mathcal{S}}\tilde{b}(y_i)| > \tfrac{1}{2}L^{-1+\eta/3}]. \end{aligned}
$$

To bound the terms of the sum on the right-hand side of (3.77), we separately estimate $\mathbb{P}[\partial_j G_{\mathcal{S}}\tilde{b}(y_i) > \frac{1}{2d}L^{-1+\eta/3}]$ and $\mathbb{P}[\partial_j G_{\mathcal{S}}\tilde{b}(y_i) < -\frac{1}{2d}L^{-1+\eta/3}]$ for $j = 1, \ldots, d$ with the help of an Azuma-type inequality. Therefore we cover the slab $\mathcal{S}$ with disjoint cubes of side-length $R$ and assign these cubes to $2^d$ disjoint families of cubes that are spaced by a distance $R$. We denote with $Q_k^m = x_{m,k} + [-\frac{R}{2}, \frac{R}{2})^d$, $1 \leq m \leq 2^d$, $k \geq 1$, the cubes associated to the $m$th family and define for $i \in \mathcal{I}, 1 \leq j \leq d, \omega \in \Omega$,

$$
(3.78) \quad Y_{i,k}^m(\omega) = \int_{Q_k^m \cap \mathcal{S}} \partial_j g(y_i, z)\tilde{b}(z, \omega)\, dz, \qquad k \geq 1.
$$

For fixed $m \in \{1, \ldots, 2^d\}$ and $i \in \mathcal{I}, 1 \leq j \leq d$, these random variables are $\mathbb{P}$-independent (as $k$ varies) and have mean 0 by Fubini's theorem. Moreover, it follows from (3.49) that for all $\omega \in \Omega$; $m, i, j, k \geq 1$,

$$
(3.79) \quad \begin{aligned} |Y_{i,k}^m(\omega)| &\leq cL^{-1}(|x_{m,k} - y_i|^{1-d} \wedge 1 + L^{1-d})\exp(-c_{17}|x_{m,k} - y_i|_\perp/L) \\ &\stackrel{\text{def}}{=} \gamma_{m,k}. \end{aligned}
$$

Indeed, either $|y_i - x_{m,k}| \leq \sqrt{d}R$ and using polar coordinates we obtain that $|Y_{i,k}^m| \leq c\varepsilon \int_{B_{2\sqrt{d}R}(y_i)}(r^{1-d} + L^{1-d})r^{d-1}\, dr \leq c'L^{-1}(1 + L^{1-d})$, or $|y_i - x_{m,k}| \geq \sqrt{d}R$ and we can bound the integral by the supremum of the integrand times

the constant volume of $Q_k^m$. Using a slight variation of the proof of Azuma's inequality, we find for $1 \leq j \leq d$, $i \in \mathcal{I}$,

$$
\begin{aligned}
\mathbb{P}\left[\partial_j G_{\mathcal{S}} \tilde{b}(y_i) > \frac{1}{2d} L^{-1+\eta/3}\right] &\leq \sum_{m=1}^{2^d} \mathbb{P}\left[\sum_{k \geq 1} Y_k^m(\omega) > \frac{1}{d2^{d+1}} L^{-1+\eta/3}\right] \\
&\leq \sum_{m=1}^{2^d} \exp\left(-\frac{d^{-2} 2^{-2(d+1)} L^{-2+2/3\eta}}{\sum_{k \geq 1} (\gamma_{m,k})^2}\right) \\
&\leq 2^d \exp(-cL^{2/3\eta}),
\end{aligned}
\tag{3.80}
$$

since the following easy computation and (3.79) show that $\sum_{k \geq 1}(\gamma_{m,k})^2$ is of order $L^{-2}$ for all $m \geq 1$:

$$
\begin{aligned}
L^2 \sum_{k \geq 1}(\gamma_{m,k})^2 &\leq c \sum_{|x_{m,k}-y_i| \leq 4L} (|x_{m,k}-y_i|^{-d+1} \wedge 1 + L^{-d+1})^2 \\
&\quad + \sum_{|x_{m,k}-y_i|_\perp \geq 2L} L^{-2d+2} \exp(-c'|x_{m,k}-y_i|_\perp/L) \\
&\leq c \int_1^{4L} (r^{-2d+2} + L^{-2d+2}) r^{d-1}\, dr \\
&\quad + L^{-2d+3} \int_L^\infty e^{-c'r/L} r^{d-2}\, dr \\
&\leq c + L^{-d+2} \int_1^\infty e^{-c'u} u\, du \leq c''.
\end{aligned}
\tag{3.81}
$$

The same bound as in (3.80) holds for the terms $\mathbb{P}[\partial_j G_{\mathcal{S}} \tilde{b}(y_i) < -\frac{1}{2d} L^{-1+\eta/3}]$, $1 \leq j \leq d$. Collecting the estimates (3.77), (3.80) and recalling that the cardinality of $\mathcal{I}$ is polynomial in $L$, we have proved the claim (3.70).

Finally we come to (3.71). The argument is similar to the previous one. First we handle the infimum over $\mathcal{V}$ by covering $\mathcal{V}$ with disjoint cubes of the form $x_i + [-\frac{R}{2}, \frac{R}{2}]^d$, for some adequate points $x_i, i \in \mathcal{I}'$ where $x_i \cdot e_1 = 0$ and $|\mathcal{I}'| \leq cL^{12(d-1)}$. Then it follows from (3.34) that for all $\omega \in \Omega$ and $|x - x_i| \leq \frac{R}{2}$,

$$
|G_{\mathcal{S}} \tilde{b}(x) - G_{\mathcal{S}} \tilde{b}(x_i)| \leq c_{15} 2\varepsilon L \frac{R}{2} \sqrt{d} \leq c.
\tag{3.82}
$$

Hence the discretization step implies that the left-hand side of (3.71) is less than

$$
\sum_{x_i} \mathbb{P}\left[G_{\mathcal{S}} \tilde{b}(x_i) \leq -\frac{c_{12}}{8} L^\eta\right].
\tag{3.83}
$$

Then we use the same $2^d$ $R$-disjoint families of boxes $Q_k^m$ as before to cover the slab $\mathcal{S}$ and we define for $i \in \mathcal{I}', m \geq 1$ and all $\omega \in \Omega$,

$$\tilde{Y}_{i,k}^m(\omega) = \int_{Q_k^m \cap \mathcal{S}} g(x_i, z) \tilde{b}(z, \omega) \, dz, \qquad k \geq 1. \tag{3.84}$$

Again we observe that for fixed $m \in \{1, \ldots, 2^d\}$ and $i \in \mathcal{I}'$, these random variables are $\mathbb{P}$-independent and have mean 0. Moreover, it follows from (3.48) that for all $\omega \in \Omega; m, i, k \geq 1$

$$|\tilde{Y}_{i,k}^m(\omega)| \leq cL^{-1}(|x_{m,k} - x_i|^{2-d} \wedge 1) \exp(-c_{17}|x_{m,k} - x_i|_\perp/L) \stackrel{\mathrm{def}}{=} \tilde{\gamma}_{m,k}. \tag{3.85}$$

A computation as in (3.81) shows that for large enough $L$ and for all $1 \leq m \leq 2^d$:

$$\sum_{k \geq 1} (\tilde{\gamma}_{m,k})^2 \leq L^{-2} \begin{cases} c \log L, & d = 4, \\ c, & d \geq 5. \end{cases} \tag{3.86}$$

Then the same Azuma-type argument as before yields for large enough $L$ that each term in (3.83) is less than

$$\begin{aligned} &\exp(-cL^{2+2\eta}/\log(L)) && \text{when } d = 4, \quad \text{respectively} \\ &\exp(-cL^{2+2\eta}) && \text{when } d \geq 5. \end{aligned} \tag{3.87}$$

This completes the proof of (3.71) and thus of Proposition 3.10. $\square$

## APPENDIX

**A.1. Proof of Lemma 2.3.** We now give the proof of Lemma 2.3. In order to bound $\rho_1(\omega)$ on $\mathcal{G}$ [see (2.16)], we first construct a function which—after appropriate normalization—dominates $P_{x,\omega}[\tilde{T}_{-L_1+R+1} < \tilde{T} \wedge T_{L_1+1}]$. For the construction, we divide the box $B_1$ into slabs of width $L_0$ and consider an expression inspired from the solution of a discrete one-dimensional Dirichlet problem for the exit probability of a Markov chain whose states correspond in essence to the slabs $\mathcal{S}_\rangle, i \in \mathbb{Z}$.

Indeed, we recall (2.13) and for integers $a < b$, we consider the products $\prod_{a,b} = \prod_{j=a+1}^b \hat{\rho}(j, \omega)^{-1}$ and set $\prod_{a,a} = 1$. Then we define the function $f$ on $\{-n_0 + 1, -n_0 + 2, \ldots, n_0 + 2\} \times \Omega$ via

$$\begin{aligned} &f(n_0 + 2, \omega) = 0, \qquad f(n_0 + 1, \omega) = 1, \\ &f(i, \omega) = \sum_{i \leq m \leq n_0 + 1} \prod_{m, n_0 + 1} \qquad \text{for } i \leq n_0. \end{aligned} \tag{A.1}$$

For simplicity we drop the $\omega$-dependence from the notation. We now show that for $\omega \in \Omega$,

$$P_{0,\omega}[\tilde{T}_{-L_1+R+1} < \tilde{T} \wedge T_{L_1+1}] \leq \frac{f(0)}{f(1 - n_0)}. \tag{A.2}$$

Let us introduce the $(\mathcal{F}_{V_m})_{m\geq 0}$ -stopping time

(A.3) $$\tau = \inf\{m \geq 0 : X_{V_m} \in \mathcal{S}_{n_0+2} \cup \mathcal{S}_{1-n_0}\}.$$

Observe that $P_{0,\omega}$-a.s. on the event which appears in (A.2), $X_{V_\tau} \in \mathcal{S}_{1-n_0}$ and $V_\tau < \tilde{T}$, and thus for $\omega \in \Omega$,

(A.4) $$P_{0,\omega}[\tilde{T}_{-L_1+R+1} < \tilde{T} \wedge T_{L_1+1}] \leq \frac{E_{0,\omega}[f(I(X_{V_\tau})), V_\tau < \tilde{T}]}{f(1-n_0)}.$$

As we will see now, the numerator on the right-hand side is less than $f(0)$: for $\omega \in \Omega$, $m \geq 0$,

(A.5) $$\begin{aligned}E_{0,\omega}[f(I(X_{V_{(m+1)\wedge\tau}})), V_{(m+1)\wedge\tau} \leq \tilde{T}] \\ \leq E_{0,\omega}[f(I(X_{V_{m\wedge\tau}})), V_{m\wedge\tau} \leq \tilde{T}, \tau \leq m] \\ + E_{0,\omega}[f(I(X_{V_{m+1}})), V_m \leq \tilde{T}, \tau > m]\end{aligned}$$

and by the strong Markov property, the second term on the right-hand side equals

(A.6) $$E_{0,\omega}[V_m \leq \tilde{T}, \tau > m, E_{X_{V_m},\omega}[f(I(X_{V_1}))]].$$

However on $\{V_m \leq \tilde{T}, \tau > m\}$, $P_{0,\omega}$-a.s.:

(A.7) $$\begin{aligned}E_{X_{V_m},\omega}&[f(I(X_{V_1}))] \\ &= f(I(X_{V_m})) + \hat{p}(X_{V_m})[f(I(X_{V_m})+1) - f(I(X_{V_m}))] \\ &\quad + \hat{q}(X_{V_m})[f(I(X_{V_m})-1) - f(I(X_{V_m}))] \\ &\stackrel{(A.1)}{=} f(I(X_{V_m})) + \prod_{I(X_{V_m}),n_0+1}[-\hat{p}(X_{V_m}) + \hat{q}(X_{V_m})\rho(I(X_{V_m}))^{-1}].\end{aligned}$$

Note that $P_{0,\omega}$-a.s., $X_{V_m} \in \mathcal{S}_{I(X_{V_m})}$, for $m \geq 0$. Hence the expression inside the square brackets is nonpositive; see (2.13). As a result, we obtain that the left-hand side of (A.5) is smaller than or equal to $E_{0,\omega}[f(I(X_{V_{m\wedge\tau}})), V_{m\wedge\tau} \leq \tilde{T}]$. The latter expression is hence nonincreasing with $m$. Since $\tau$ is $P_{0,\omega}$-a.s. finite, it follows from Fatou's inequality that for $\omega \in \Omega$,

(A.8) $$E_{0,\omega}[f(I(X_{V_\tau})), V_\tau \leq \tilde{T}] \leq f(0).$$

Together with (A.4), this implies (A.2).

We now derive a bound on $\rho_1$. Let us define for $\omega \in \Omega$,

(A.9) $$A = P_{0,\omega}[\tilde{T}_{-L_1+R+1} < \tilde{T} \wedge T_{L_1+1}] + P_{0,\omega}[\tilde{T} < \tilde{T}_{-L_1+R+1} \wedge T_{L_1+1}].$$

Observe that $q(0,\omega) \leq A$ and since $\frac{q}{1-q}$ is nondecreasing in $q$, we obtain for $\omega \in \Omega$ that $\rho_1(\omega) \leq \frac{A}{(1-A)_+}$. Using (A.2) and (2.16), it follows for $\omega \in \mathcal{G}$ that

(A.10) $$\rho_1(\omega) \leq \frac{f(0) + f(1-n_0)\kappa^{9L_1}}{(f(1-n_0) - f(0) - f(1-n_0)\kappa^{9L_1})_+}.$$

Let us for the time being assume that there is a $c_1 > R + 2$ such that for $L_0 \geq c_1$ and $\omega \in \Omega$,

(A.11) $$f(0) + f(1 - n_0)\kappa^{9L_1} \leq 2f(0),$$

(A.12) $$f(1 - n_0) - f(0) - f(1 - n_0)\kappa^{9L_1} \geq \prod_{-n_0+1, n_0+1}.$$

Then in view of (A.10) and the definition of $f(0)$, for $L_0 \geq c_1$, $\omega \in \mathcal{G}$,

(A.13) $$\rho_1(\omega) \leq 2 \sum_{0 \leq m \leq n_0+1} \prod_{-n_0+1 < j \leq m} \hat{\rho}(j, \omega).$$

Observe that by the definition (2.13), $\{\hat{\rho}(j,\omega), j \text{ even}\}$ and $\{\hat{\rho}(j,\omega), j \text{ odd}\}$ are two collections of independent random variables, as $\rho(j, \omega)$ and $\rho(j+2, \omega)$ depend on regions separated by a distance $R$. With the help of Cauchy–Schwarz's inequality and $(u + v)^a \leq u^a + v^a$, for $u, v \geq 0$ and $a \in (0, 1]$, we find that for $L_0 \geq c_1$,

(A.14) $$\mathbb{E}[\rho_1^a, \mathcal{G}] \leq 2 \sum_{0 \leq m \leq n_0+1} \prod_{-n_0+1 < j \leq m} \mathbb{E}[\hat{\rho}(j, \omega)^{2a}]^{1/2}.$$

From Lemma 2.1, we have that for all $\omega \in \Omega$, $\rho_1(\omega) \leq \kappa^{-L_1-1}$. This inequality and (A.14) immediately imply the claim (2.15).

Let us now show (A.11). Using again Lemma 2.1, we have for all $\omega \in \Omega$, $-n_0 + 1 \leq j \leq n_0 + 1$:

(A.15) $$\kappa^{L_0+1} \leq \hat{\rho}(j, \omega) \leq \kappa^{-(L_0+1)}.$$

In view of (A.1) and since $L_0 + 1 \leq 2L_0$, we find that

$$f(1 - n_0)\kappa^{9L_1} \leq (2n_0 + 1)\kappa^{-(L_0+1)2n_0}\kappa^{9L_1} \leq (2n_0 + 1)\kappa^{5n_0 L_0}.$$

If $L_0 \geq c_1 \geq R + 2$ large enough, it follows that for all $\omega \in \Omega$ and all $n_0 \geq 3$,

(A.16) $$f(1 - n_0)\kappa^{9L_1} \leq \kappa^{4n_0 L_0} < 1.$$

Clearly $f(0) \geq 1$ and we obtain (A.11). To see (A.12), we note that

(A.17) $$\begin{aligned} f(1 - n_0) - f(0) &\geq \prod_{-n_0+1, n_0+1} + \prod_{-1, n_0+1} \\ &\overset{(A.15)}{\geq} \prod_{-n_0+1, n_0+1} + \kappa^{(L_0+1)(n_0+2)}. \end{aligned}$$

Since $(L_0 + 1)(n_0 + 2) \leq 4L_0 n_0$ and because of (A.16), the claim (A.12) follows, provided that $L_0 \geq c_1$. This finishes the proof of Lemma 2.3.

**A.2. Proof of Lemma 3.9.** We now prove Lemma 3.9. We start with the proof of (3.49). A similar and easier argument also shows (3.48). Since for $d \geq 4$, we have

$$(A.18) \quad |\partial_i g_d(x,y)| \leq c|x-y|^{1-d} \quad \text{and} \quad |\partial_i \partial_j g_d(x,y)| \leq c'|x-y|^{-d},$$

the sum of the first and second derivatives of the terms with $k \geq 2$ appearing in (3.30) converges uniformly for all $x, y \in \mathcal{S}$. Hence $g(x,y)$ is twice continuously differentiable for $x, y \in \mathcal{S}$, $x \neq y$, and interchanging differentiation and summation yields for all $x, y \in \mathcal{S}$

$$(A.19) \quad |\nabla g(x,y)| \leq 3|\nabla g_d(x,y)| + c(2L)^{-d+1} \quad \text{and}$$

$$(A.20) \quad |\partial_i \partial_j g(x,y)| \leq 3|\partial_i \partial_j g_d(x,y)| + c'(2L)^{-d}, \quad \text{as well as}$$

$$(A.21) \quad \Delta g(x,y) = 0 \qquad \text{for } x \neq y.$$

For any $x \in \mathcal{S}$, we consider a small vector $h$ with $x+h \in \mathcal{S}$ and an point $y \in \mathcal{S}$ with $|x-y|_\perp \geq L$. Moreover, we denote with $W$ a $d$-dimensional Brownian motion starting at $y$ under some measure $P$ and with $T$ the stopping time $\inf\{t \geq 0; |W_t - y|_\perp \geq \frac{1}{2}|x-y|\}$. Since $g(x,y)$ is symmetric in $x$ and $y$, it is also harmonic in $y$ and thus $g(x, W_{t \wedge T \wedge T_\mathcal{S}})$ is a bounded martingale under $P$. The stopping theorem thus implies that

$$(A.22) \quad \begin{aligned} &\frac{1}{|h|}|g(x+h,y) - g(x,y)| \\ &= \frac{1}{|h|}|E_P[g(x+h, W_{T \wedge T_\mathcal{S}}) - g(x, W_{T \wedge T_\mathcal{S}})]|. \end{aligned}$$

Direct inspection of $g(x,y)$ shows that it vanishes on the boundary of $\mathcal{S}$. Hence using the mean value theorem, the latter expression is smaller than

$$(A.23) \quad \sup\{|\nabla g(x',y')|; x' \in B(x,h), |y'-y|_\perp = \tfrac{1}{2}|x-y|\} P[T < T_\mathcal{S}].$$

Because of (A.19) the first factor above is less than $c|x-y|^{1-d} + c'L^{1-d}$ and a scaling argument similar to the one leading to (3.60) yields that $P[T < T_\mathcal{S}] \leq \exp(-c\frac{|x-y|_\perp}{L})$. Letting $h$ tend to 0 in (A.22), (A.23) and treating the cases $|x-y| < L$ and $|x-y| > L$ separately, we obtain the claimed result (3.49). The same martingale argument also leads to (3.48).

We now prove (3.50). For any $x_0 \in \mathcal{S}$, we define the auxiliary set $U = \{x \in \mathcal{S}; |x-x_0|_\perp < 1\}$. From (3.31) we can write the Green operator for a bounded Hölder continuous function $f$ as follows: we define $\tilde{g} = g - g_d$ and for any $x \in U$, we find

$$(A.24) \quad \begin{aligned} G_\mathcal{S} f(x) &= \int_U g_d(x,y) f(y)\, dy + \int_U \tilde{g}(x,y) f(y)\, dy \\ &\quad + \int_{\mathcal{S} \setminus U} g(x,y) f(y)\, dy. \end{aligned}$$

According to [5], Lemma 4.2, the first term on the right-hand side is twice continuously differentiable on $U$, and its Laplacian equals $-2f(x)$. With the same argument as below (A.18), we see that $\tilde{g}(\cdot, y)$ is harmonic on $U$ for any $y \in \mathcal{S}$. Hence Fubini's theorem together with the mean value theorem (see [5], Theorem 2.7) yield that the second term on the right-hand side of (A.24) is harmonic on $U$. The same is valid for the last term, since from (A.21), $g(\cdot, y)$ is harmonic on $U$ for any $y \in \mathcal{S} \setminus U$. As $x_0 \in \mathcal{S}$ is arbitrary, we obtain (3.50). This finishes the proof of Lemma 3.9.

**Acknowledgments.** I would like to thank Professor A.-S. Sznitman for guiding me through this work with patience and constant advice. I am also grateful for many helpful discussions with Tom Schmitz.

## REFERENCES

[1] ALON, N. and SPENCER, J. H. (1992). *The Probabilistic Method*. Wiley, New York. MR1140703
[2] BASS, R. F. (1997). *Diffusions and Elliptic Operators*. Springer, New York. MR1483890
[3] BENSOUSSAN, A., LIONS, J. L. and PAPANICOLAOU, G. (1978). *Asymptotic Analysis for Periodic Structures*. North-Holland, Amsterdam. MR0503330
[4] BOLTHAUSEN, E. and GOLDSHEID, I. (2000). Recurrence and transience of random walks in random environments on a strip. *Comm. Math. Phys.* **214** 429–447. MR1796029
[5] GILBARG, D. and TRUDINGER, N. S. (1998). *Elliptic Partial Differential Equations of the Second Order*. Springer, Berlin. MR1814364
[6] HUGHES, B. D. (1996). *Random Walks in Random Environment* **2**. Clarendon Press, Oxford. MR1420619
[7] KALIKOW, S. A. (1981). Generalized random walk in a random environment. *Ann. Probab.* **9** 753–768. MR0628871
[8] KARATZAS, I. and SHREVE, S. E. (1991). *Brownian Motion and Stochastic Calculus*. Springer, Berlin. MR1121940
[9] MOLCHANOV, S. (1994). Lectures on random media. *Lectures on Probability Theory (Saint-Flour, 1992)*. *Lecture Notes in Math.* **1581** 242–411. Springer, Berlin. MR1307415
[10] REVESZ, P. (1990). *Random Walk in Random and Non-Random Environments*. World Scientific, Singapore. MR1082348
[11] SCHMITZ, T. (2006). Diffusions in random environment and ballistic behavior. *Ann. Inst. H. Poincaré Probab. Statist.* **42** 683–714. MR2269234
[12] SCHMITZ, T. (2006). Examples of condition $(T)$ for diffusions in random environment. *Electron. J. Probab.* **11** 540–562. MR2242655
[13] SHEN, L. (2003). On ballistic random diffusions in random environment. *Ann. Inst. H. Poincaré Probab. Statist.* **39** 839–876. MR1997215
[14] SHEN, L. (2004). Addendum to the article "On ballistic random diffusions in random environment." *Ann. Inst. H. Poincaré Probab. Statist.* **40** 385–386. MR2060459
[15] SOLOMON, F. (1975). Random walks in a random environment. *Ann. Probab.* **3** 1–31. MR0362503
[16] STROOCK, D. W. (1993). *Probability Theory*; *An Analytic View*. Cambridge Univ. Press. MR1267569

[17] STROOCK, D. W. and VARADHAN, S. R. S. (1979). *Multidimensional Diffusion Processes*. Springer, Berlin. MR0532498
[18] SZNITMAN, A.-S. (2001). A class of transient random walks in random environment. *Ann. Probab.* **29** 723–764. MR1849176
[19] SZNITMAN, A.-S. (2002). An effective criterion for ballistic behavior of random walks in random environments. *Probab. Theory Related Fields* **122** 509–544. MR1902189
[20] SZNITMAN, A.-S. (2003). On new examples of ballistic random walks in random environment. *Ann. Probab.* **31** 285–322. MR1959794
[21] SZNITMAN, A.-S. and ZERNER, M. (1999). A law of large numbers for random walks in random environment. *Ann. Probab.* **27** 1851–1869. MR1742891
[22] ZEITOUNI, O. (2004). Random walks in random environment. *Lectures on Probability Theory and Statistics. Lecture Notes in Math.* **1837** 189–312. Springer, Berlin. MR2071631

DEPARTMENT OF MATHEMATICS
ETH ZURICH
CH-8092 ZURICH
SWITZERLAND
E-MAIL: goergen@math.ethz.ch